\documentclass{elsarticle}

\usepackage{amsmath}
\usepackage{amssymb}
\usepackage{amsthm}
\usepackage{accents}
\usepackage{bbm}
\usepackage{bm} 
\usepackage{lineno,hyperref}
\usepackage{color,soul}  
\usepackage{flafter}
\usepackage{fullpage}
\usepackage{graphicx}
\usepackage{epstopdf}
\usepackage{latexsym}
\usepackage{subfigure}
\usepackage{enumerate}
\usepackage[mathscr]{eucal}
\usepackage[small]{caption}

\newif\ifupdatetikz
\updatetikzfalse 
\ifupdatetikz
	\usepackage{tikz}
	\usepackage{pgfplots}
	\pgfplotsset{
		compat=newest,
		tick label style={font=\scriptsize},
		label style={font=\scriptsize},
		legend style={font=\scriptsize}
	}
	\usepgfplotslibrary{external}
	\usetikzlibrary{shapes,arrows,positioning,calc}
	\tikzset{cross/.style={cross out, draw=black, minimum size=2*(#1-\pgflinewidth), inner sep=0pt, outer sep=0pt},cross/.default={2pt}}
	\tikzstyle{red_dot} = [red,fill=red]
	\tikzstyle{axis} = [-latex',line width=1.25]
	
\else
	\usepackage{tikzexternal}
\fi
\newcommand{\figname}[1]{\tikzsetnextfilename{#1}}

\newtheoremstyle{mystyle}%
  {}
  {}
  {}
  {}
  {\bfseries}
  {.}
  { }
  {\thmname{#1}\thmnumber{ #2}\thmnote{ (#3)}}

\theoremstyle{mystyle}

\newtheorem{thm}{Theorem}[section]

\newtheorem{rem}[thm]{Remark}

\let\originalleft\left
\let\originalright\right
\renewcommand{\left}{\mathopen{}\mathclose\bgroup\originalleft}
\renewcommand{\right}{\aftergroup\egroup\originalright}

\newcommand{\Dt}{\Delta t}

\newcommand{\dt}{\delta t}

\newcommand{\E}{\mathbb{E}}					
\newcommand{\epsi}{\varepsilon}

\newcommand{\Fbar}{\bar{F}}
\newcommand{\Fhat}{\widehat{F}}
\newcommand{\FHMM}[3]{\Fhat^{\HMM}_{#2}(#1;#3)}			
\newcommand{\FHMMbar}[2]{\Fbar^{\HMM}_{#2}(#1)}			

\newcommand{\HMM}{\text{\scriptsize{HMM}}}		

\newcommand{\lipschitzCst}{L_c}

\newcommand{\minv}{m_\infty}
\newcommand{\muinv}{\mu^\infty}

\newcommand{\normal}[2]{\mathcal{N}\left(#1,#2\right)}	

\newcommand{\R}{\mathbb{R}}				

\newcommand{\rhoinv}{\rho^\infty}

\newcommand{\set}[2]{\left(#1\right)_{#2=1}^{\uppercase{#2}}}
\newcommand{\siginv}{\sigma^2_\infty}

\newcommand{\Var}{\mathrm{Var}}

\newcommand{\Xbar}{\bar{X}}
\newcommand{\Xhat}{\widehat{X}}
\newcommand{\Y}{\mathcal{Y}}

\newlength{\dhatheight}

\newcommand{\Atilde}{\widetilde{A}}
\newcommand{\Btilde}{\widetilde{B}}

\begin{document}

\begin{frontmatter}

\title{Variance-reduced multiscale simulation of slow-fast stochastic differential equations}

\author[mainAddress]{Ward Melis\corref{correspondingAuthor}}
\cortext[correspondingAuthor]{Corresponding author}
\ead{ward.melis@cs.kuleuven.be}

\author[mainAddress]{Giovanni Samaey}
\ead{giovanni.samaey@cs.kuleuven.be}

\address[mainAddress]{Department of Computer Science, K.U. Leuven, Celestijnenlaan 200A, 3001 Leuven, Belgium}

\begin{abstract}
	We study a variance reduction strategy based on control variables for simulating the averaged macroscopic behavior of a stochastic slow-fast system. We assume that this averaged behavior can be written in terms of a few slow degrees of freedom, and that the fast dynamics is ergodic for every fixed value of the slow variable. The time derivative for the averaged dynamics can then be approximated by a Markov chain Monte Carlo method. The variance-reduced scheme that is introduced here uses the previous time instant as a control variable. We analyze the variance and bias of the proposed estimator and illustrate its performance when applied to a linear and nonlinear model problem.
\end{abstract}


\end{frontmatter}


\section{Introduction} \label{sec:intro}
Stochastic differential equations (SDEs) are ubiquitous in a multitude of real-life applications, appearing in different scientific domains such as climate and environmental sciences \cite{Berglund2006,Bruna2014,Imkeller2001}, molecular dynamics \cite{Erban2009,Givon2004} and bacterial chemotaxis \cite{Rousset2013}. Many of these applications contain processes that inherently evolve over multiple time scales, leading to excessive computational cost with standard time discretization methods. As a consequence, there is currently a large interest in developing dedicated numerical methods that circumvent, or even exploit, the presence of a time-scale separation in the problem at hand.

Clearly, the development of new numerical techniques needs to be supplemented by a detailed analysis of their efficiency and accuracy, and this for a set of model problems of which the multiscale nature is well understood. One prototypical example system that was proposed in \cite{Givon2004} to analyze such convergence behavior is a singularly perturbed slow-fast system in which the slow variable is described deterministically, while the model for the fast variable contains stochastic effects. The specific form is as follows:
\begin{equation} \label{eq:slow_fast_system}
	\begin{cases} 
	dx(t) = f(x,y)dt, &x(0) = x_0 \in \R \\
	dy(t) = \dfrac{1}{\epsi}g(x,y)dt + \dfrac{1}{\sqrt{\epsi}}\beta(x,y)dW(t), &y(0) = y_0 \in \R.
	\end{cases}
\end{equation}
where the scalar quantities $x(t):[0,T] \to \R$ and $y(t):[0,T] \to \R$ represent the slow and fast evolving stochastic processes, respectively. The functions $f(x,y),g(x,y) \in \R$ are called the drift functions and $\beta(x,y) \in \R$ is termed the diffusion function. Furthermore, $W(t) \in \R$ denotes a standard Brownian motion. 
The parameter $\epsi \ll 1$ is a positive small-scale parameter that measures the time scale separation between the fast and slow variable in system \eqref{eq:slow_fast_system}. 
In addition, we assume that the fast dynamics is ergodic for every fixed state $X \in \R$ of the slow variable, implying the existence and uniqueness of an invariant measure \cite{Pavliotis2008}. 
We note that the differential form used in system \eqref{eq:slow_fast_system} is purely formal, given that Brownian paths are continuous everywhere but nowhere differentiable. Consequently, system \eqref{eq:slow_fast_system} should be understood in the integral form, where stochastic integrals are interpreted in the It\^o-sense. In general, the SDE may be very high-dimensional(especially with many fast degrees of freedom), see, for instance, \cite{Givon2004}.

Often, one is only interested in the evolution of the slow variable of system \eqref{eq:slow_fast_system} and not in the detailed evolution of the fast variable. However, the fast dynamics cannot be omitted, since the slow process explicitly depends on the fast variable. 
Due to the stiffness in system \eqref{eq:slow_fast_system}, explicit simulation techniques such as the Euler-Maruyama or higher-order Milstein schemes are computationally prohibitive. 
Also, implicit methods fail to capture the correct invariant measure, thus introducing a bias, see \cite{Li2008}.

The difficulties related to direct simulation can be avoided by exploiting the time scale separation of system \eqref{eq:slow_fast_system}: for $\epsi \to 0$, the averaging principle yields a reduced description for the slow variable:
\begin{equation} \label{eq:averaged_eq} 
	\frac{dX}{dt} = F(X), \qquad F(X) = \int_\Y f(X,y) d\muinv_{X}(y),
\end{equation}
in which $\muinv_{X}(y)$ denotes the invariant measure induced by the fast dynamics of system \eqref{eq:slow_fast_system} keeping $x=X$ fixed, see, for instance, \cite{Pavliotis2008} and references therein. Equation \eqref{eq:averaged_eq} is known as the \textit{averaged}, \textit{macroscopic} or \textit{reduced} equation for the slow variable.

Based on this averaged equation, a method for \eqref{eq:slow_fast_system} was proposed in \cite{Vanden-Eijnden2003} and analyzed in \cite{E2005}. It consists of a macroscopic solver, such as the forward Euler or a higher-order Runge-Kutta method, to simulate \eqref{eq:averaged_eq}, combined with a procedure to estimate the effective force $F(X)$ in equation \eqref{eq:averaged_eq}. If the invariant measure is known explicitly and can readily be sampled, the integral in \eqref{eq:averaged_eq} can be approximated by a direct Monte Carlo estimator (see, e.g., \cite{caflisch1998monte} and references therein).
In general, however, the invariant measure is not known explicitly. Then, one may resort to a Markov chain Monte Carlo method, as is done in \cite{Vanden-Eijnden2003}. 
This method fits in the class of \emph{heterogeneous multiscale methods} (HMM) that were introduced in \cite{E2003a} for a broad class of multiscale problems and provide a natural setting for numerical analysis, see also \cite{Abdulle2012} for a recent review. Similar methods have been introduced based on the concept of \emph{coarse projective integration} \cite{Givon2006}. There, instead of performing one (or a few) long Markov chain Monte Carlo simulation of the fast equation, one initializes a \emph{large ensemble} of realizations, which are simulated on a short time interval. This method falls in the class of \emph{equation-free methods} \cite{Kevrekidis2003,Kevrekidis2009}, see also \cite{Gear2002}, and can also be used in a more general setting where one is unable to identify or constrain the slow degree of freedom.

Unfortunately, the statistical error of the above-described methods can be quite large, and decreases only as $M^{-1/2}$ when the number of samples $M$ tends to infinity.
In this work, which expands the results reported in \cite{melis2016variance}, we therefore propose a variance-reduction technique based on control variables, see, e.g., \cite{caflisch1998monte,Glasserman2003}. The method can be applied both in the projective integration and the HMM setting, and bears some resemblance to the technique that was proposed in \cite{Papavasiliou2007} for variance-reduced coarse projective integration of SDEs of the form \eqref{eq:slow_fast_system}. 
While we present the method and main analysis in the HMM framework, we will comment on coarse projective integration where appropriate.
The control variable that we introduce is based on correlating estimations of the time derivative $F(X)$ in equation~\eqref{eq:averaged_eq} on different time instants. 

The remainder of this paper is structured as follows. In section \ref{sec:slow_fast}, we introduce the stochastic slow-fast system that we intend to solve numerically. In section \ref{sec:numerical}, we describe the HMM framework to efficiently integrate these slow-fast systems. In that section, we also present the variance-reduced HMM method, and we comment on the applicability of this method for coarse projective integration. Next, in section \ref{sec:num_prop} we analyze the numerical properties of the proposed variance reduction method. Numerical results are reported in section \ref{sec:results}. We conclude in section \ref{sec:conclusions} with a brief discussion and ideas for future work.


\section{Slow-fast system} \label{sec:slow_fast}
The general form of slow-fast systems we consider in this work is given in equation \eqref{eq:slow_fast_system}. In what follows, we will always assume that the fast dynamics of system \eqref{eq:slow_fast_system} is ergodic for all fixed values of the slow variable. This means that the fast equation produces a unique invariant measure for every fixed value $X$, denoted by $\muinv_{X}(y)$. Ergodicity implies that the statistical properties of the ensemble of the stochastic process at a fixed time instant and those of one realization of the process over an infinite time interval are the same. Consequently, for an ergodic process, averaging a function with respect to the invariant measure yields the same result as averaging this function over one infinitely long time path of the process:
\begin{equation} \label{eq:ergodicity_property}
	F(X) = \int_\Y f(X,y)d\muinv_X(y) = \lim_{T \to \infty}\frac{1}{T}\int_0^T f(X,y(t+\tau))d\tau.
\end{equation}
Equation \eqref{eq:ergodicity_property} serves as a base for numerical methods avoiding explicit knowledge of the invariant measure $\muinv_{X}(y)$.
Additionally, we will always assume that the invariant measure possesses a density $\rhoinv_{X}$ with respect to the Lebesgue measure: $d\muinv_{X}(y) = \rhoinv_{X}(y)dy$. 

Moreover, we assume that the function $F(X)$ in equation \eqref{eq:averaged_eq} is Lipschitz continuous with Lipschitz constant $\lipschitzCst$, implying the following inequality:
\begin{equation} \label{eq:F_lipschitz}
	|F(X_1) - F(X_2)| \le \lipschitzCst|X_1 - X_2|, \qquad \forall X_1,X_2 \in \R,
\end{equation}
and we assume the functions $f$, $g$ and $\beta$ sufficiently differentiable such that all derivatives exist that are required during the analysis.

In the following two paragraphs, we introduce the linear and nonlinear illustrative examples that will be used in the numerical experiments throughout the text.

\paragraph{Linear system}\label{sec:slow_fast_lin}
In the linear setting, system \eqref{eq:slow_fast_system} takes on the following form:
\begin{equation} \label{eq:lin_slow_fast}
	\begin{cases} 
	dx(t) = \big(\lambda x(t) + py(t)\big)dt  \\
	dy(t) = \dfrac{1}{\epsi}\big(qx(t) - Ay(t)\big)dt + \dfrac{1}{\sqrt{\epsi}}dW(t),
	\end{cases}
\end{equation}
in which the parameters $\lambda$, $p$, $q$ and $A$ are all real scalars. In addition, to ensure that solutions decay exponentially with time, we require that $\lambda < 0$ and $A\in\left(\dfrac{pq}{-\lambda},2\right]$.

For this linear system, the fast equation corresponds to a linear Ornstein-Uhlenbeck process with parameters $q$ and $A$ for which the invariant measure 
can be calculated analytically as \cite{Pavliotis2008}:
\begin{equation} \label{eq:lin_inv_measure} 
	\muinv_X(y) \sim \normal{\minv}{\siginv}, \qquad \minv = \frac{q}{A}X, \qquad \siginv = \frac{1}{2A}, 
\end{equation}
where $\normal{\cdot}{\cdot}$ represents the normal distribution and $\minv$ and $\siginv$ denote the mean and variance of the invariant measure, respectively. Since the invariant measure is known, the integral in equation \eqref{eq:averaged_eq} can be calculated analytically, yielding:
\begin{align} \label{eq:lin_F_exact} 
	F(X) &= \int_\Y(\lambda X + py)d\muinv_X(y) = \lambda X + p\int_\Y yd\muinv_X(y)
	= \left(\lambda + \frac{pq}{A}\right)X.
\end{align}
In that case, the reduced equation \eqref{eq:averaged_eq} for system \eqref{eq:lin_slow_fast} becomes:
\begin{equation} \label{eq:lin_averaged_eq} 
	\frac{dX}{dt} = \left(\lambda + \frac{pq}{A}\right)X,
\end{equation}
which is a linear ODE. The exact solution to equation \eqref{eq:lin_averaged_eq} with initial condition $X(0) = x_0$ is then given as:
\begin{equation} \label{eq:lin_X_exact} 
	X(t) = x_0\exp\left(\left(\lambda + \frac{pq}{A}\right)t\right).
\end{equation}

\paragraph{Nonlinear system}
As a second example, we consider the following nonlinear stochastic multiscale system from \cite{Papavasiliou2007}:
\begin{equation} \label{eq:nonlin_slow_fast}
	\begin{cases} 
	dx(t) = -\big(y(t) + y(t)^2\big)dt  \\
	dy(t) = -\dfrac{1}{\epsi}\big(y(t) - x(t)\big)dt + \dfrac{1}{\sqrt{\epsi}}dW(t).
	\end{cases}
\end{equation}
In this case, the dynamics of the slow variable is nonlinear, while the fast variable is again described by a linear Ornstein-Uhlenbeck process. Using equation \eqref{eq:lin_inv_measure}, we obtain a Gaussian invariant measure $\muinv_X(y)$ with the invariant mean $\minv=X$ and variance $\siginv=1/2$.
Since the invariant measure is known explicitly, the expression of $F$ in equation \eqref{eq:averaged_eq} can be calculated analytically as:
\begin{equation} \label{eq:nonlin_F_exact}
	F(X) = -\int_\Y yd\muinv_X(y) -\int_\Y y^2d\muinv_X(y)=
	-\left(X+X^2+\frac{1}{2}\right),
\end{equation}
resulting in a nonlinear ODE. The exact solution of the resulting macroscopic equation with initial condition $X(0) = x_0$ is obtained as:
\begin{equation} \label{eq:nonlin_X_exact}
	X(t) = -\frac{1}{2} - \frac{1}{2}\tan\left(\frac{t}{2}-\arctan(2x_0+1)\right).
\end{equation}


\section{Numerical method} \label{sec:numerical}
In this section, we construct a variance-reduced numerical scheme to solve the averaged equation \eqref{eq:averaged_eq} for the slow variable of the underlying slow-fast system given in \eqref{eq:slow_fast_system}. Since the averaged equation is deterministic, any stable explicit ODE solver can be used. 
Here, we employ the forward Euler (FE) method. To that end, we discretize equation \eqref{eq:averaged_eq} on a uniform time mesh with time step $\Dt$, and $t^n=n\Dt$. The numerical solution on this mesh is denoted by $\Xhat^n$. The forward Euler scheme for \eqref{eq:averaged_eq} is then given by,
\begin{equation} \label{eq:averaged_eq_FE}
	\Xhat^{n+1} = \Xhat^n + \Dt\Fhat(\Xhat^n), \qquad \Xhat^0 = x_0.
\end{equation}
In equation \eqref{eq:averaged_eq_FE}, the function $F$ is replaced by an appropriate estimator $\Fhat$, since, in general, the integral in equation \eqref{eq:averaged_eq} can not be calculated analytically. 


We first introduce the HMM estimator in section \ref{subsec:hmm}, where we also briefly comment on its relation to coarse projective integration. Then, we present the variance-reduced HMM estimator, which forms the focus of this paper, in section \ref{subsec:hmm_cv}.

\subsection{Heterogeneous multiscale method (HMM)} \label{subsec:hmm}
The heterogeneous multiscale method \cite{Vanden-Eijnden2003} bypasses explicit knowledge of the invariant measure in equation \eqref{eq:averaged_eq} by exploiting the ergodicity property given in equation \eqref{eq:ergodicity_property}: $F$ is calculated by averaging over one infinitely long time path of the fast process of system \eqref{eq:slow_fast_system} while keeping the value of the slow variable fixed. 
As a result, the HMM estimator boils down to a Markov chain Monte Carlo estimator: the integral in equation \eqref{eq:averaged_eq} is approximated by a Monte Carlo method, in which the samples are not drawn from the (unknown) invariant measure, but are instead generated from a Markov chain. This chain is obtained by simulating the fast equation using the explicit Euler-Maruyama scheme, which is the stochastic counterpart of the forward Euler scheme \cite{Higham2001}. 
In that regard, we discretize the fast equation on a uniform time mesh with time step $\dt$. For a given fixed value $\Xhat^n$ of the slow variable, the numerical solution at time $t^{n,m'}=n\Dt+m'\dt$ is denoted by $y^{n,m'}$. The Euler-Maruyama scheme is given by,
\begin{equation} \label{eq:EM_scheme}
	y^{n,m'+1} = y^{n,m'} + \frac{\dt}{\epsi}g(\Xhat^n,y^{n,m'}) + \sqrt{\frac{\dt}{\epsi}}\beta(\Xhat^n,y^{n,m'})\xi_n^{m'}, \qquad m'=0,...,M-1,
\end{equation}
in which $(\xi_{m'}^n)_{m'=0}^{M-1}$ is a set of mutually independent samples drawn from the standard normal distribution using a random number generator with seed $\omega_n$. 
The initial condition of the Euler-Maruyama method is chosen as $y^{0,0}=y_0$, and for all other $n>0$ as $y^{n,0}=y^{n-1,M-1}$. Then, the samples generated by the Markov chain \eqref{eq:EM_scheme} are approximately distributed according to the desired invariant measure. To eliminate the time discretization error that the Euler-Maruyama scheme induces in the invariant measure, one could add a Metropolis accept/reject step, as in the MALA algorithm \cite{roberts1996exponential}.

Since it is more natural to label samples from 1 to $M$, we use the trivial substitution $m = m' + 1$ as sample index.
Then, the HMM estimator at time instant $t^n$ using $M$ samples is calculated as follows:
\begin{equation} \label{eq:F_hmm}
	\FHMM{\Xhat^n}{M}{\omega_n} = \frac{1}{M}\sum_{m=1}^M f(\Xhat^n,y^{n,m}),
\end{equation}
in which $\omega_n$ represents the seed that is used in the random number generator. Since we can only generate finite sample sizes $M$, the HMM estimator in equation \eqref{eq:F_hmm} is a random variable.

\begin{rem}[Coarse projective integration]
	The coarse projective integration (CPI) method that was presented in \cite{Givon2006} is very similar to the method above, with a different starting point. In \cite{Givon2006}, one does not assume to be able to simulate the fast equation separately.  Instead, to advance from $\Xhat^n$ to $\Xhat^{n+1}$, one only performs short-term simulations over a time interval of size, say, $K\delta t$ with the original system~\eqref{eq:slow_fast_system}, starting from \emph{an ensemble} of initial conditions $\left\{\left(\Xhat^n,y^{n,m}\right)\right\}_{m=1}^{M}$, yielding the time-evolved ensemble $\left\{\left(\Xhat^n_{K\delta t},y^{n,m}_{K\delta t}\right)\right\}_{m=1}^{M}$. The time derivative estimator can then be obtained as 
	\begin{equation} \label{eq:F_cpi}
		\Fhat^{\textrm{CPI}}_M(\Xhat^n) = \dfrac{1}{M}\sum_{m=1}^M \dfrac{\Xhat^n_{K\delta t}-\Xhat^n}{K\delta t}.
	\end{equation}
	While the HMM and CPI methods result in somewhat different equations and have different parameters that can be chosen, the schemes are very similar.  In particular, when generating the ensemble of initial conditions for CPI using the Euler-Maruyama method~\eqref{eq:EM_scheme} and choosing the number of microscopic time steps $K=1$ in the CPI method, both methods can be seen to be identical.  The results that are obtained in this paper for the HMM method can therefore easily be carried over to the CPI case.
\end{rem}

\subsection{Variance-reduced HMM} \label{subsec:hmm_cv}

The statistical error on the estimator~\eqref{eq:F_hmm} decays only slowly (as $M^{-1/2}$) with increasing sample size $M$, which is typical for any Monte Carlo-based estimator. Moreover, in the HMM method, the samples $\set{y^{n,m}}{m}$ generated by the Markov chain \eqref{eq:EM_scheme} are clearly correlated, resulting in a higher statistical error than that of a Monte Carlo estimator using independent samples. We refer to \cite{Cances:EsaimM2An:2007} for an overview on the convergence of Markov chain Monte Carlo sampling.  Whenever the statistical error dominates the systematic error, one should reduce the variance. 
Here, we propose a variance-reduced estimator based on the control variable technique, see, e.g., \cite{caflisch1998monte}. 

\paragraph{Main idea} To clearly distinguish between the standard and variance-reduced methods, we will always denote variance-reduced estimates with an overbar, whereas standard estimates will be indicated with a hat.  Thus, with the variance-reduced estimator, the forward Euler scheme~\eqref{eq:averaged_eq_FE} becomes:
\begin{equation} \label{eq:averaged_eq_FE_cv}
	\Xbar^{n+1} = \Xbar^n + \Dt\FHMMbar{\Xbar^n}{M}, \qquad \Xbar^0 = x_0.
\end{equation}
We define the variance-reduced HMM estimator at time instant $t^n$ using $M$ samples as follows:
\begin{equation} \label{eq:F_cv}
	\FHMMbar{\Xbar^n}{M} = \FHMM{\Xbar^n}{M}{\omega_n} -  \left(\FHMM{\Xbar^{n-1}}{M}{\omega_n} - \FHMMbar{\Xbar^{n-1}}{M}\right),
\end{equation}
in which an overbar denotes a variance-reduced estimator. The first term $\FHMM{\Xbar^n}{M}{\omega_n}$ in equation \eqref{eq:F_cv} coincides with the classical HMM estimator without variance reduction for the slow variable $\Xbar^n$ at the current time instant using seed $\omega_n$. The second term $\FHMM{\Xbar^{n-1}}{M}{\omega_n}$ represents another HMM estimation without variance reduction for the slow variable $\Xbar^{n-1}$ at the previous time instant. However, this term uses the same seed $\omega_n$ as the first term and therefore $\FHMM{\Xbar^{n-1}}{M}{\omega_n}$ and $\FHMM{\Xbar^n}{M}{\omega_n}$ will be strongly correlated. The variance reduction is achieved by subtracting these two terms in an attempt to cancel out the corresponding statistical variations. The last term $\FHMMbar{\Xbar^{n-1}}{M}$ is the variance-reduced HMM estimator calculated during the previous time step, which needs to be added to avoid introducing a bias.

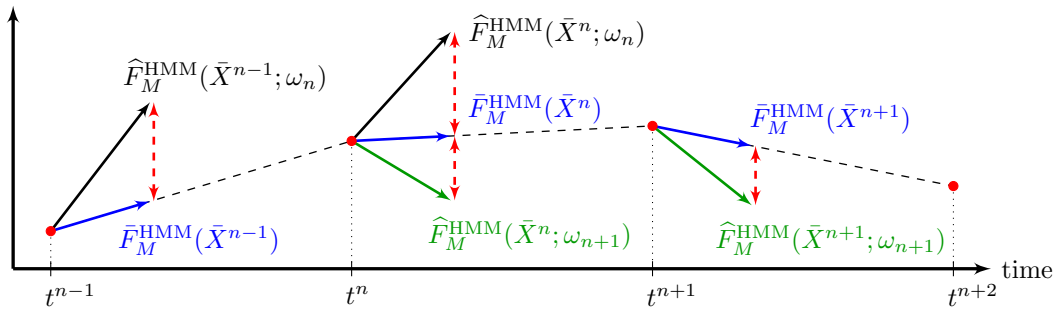
\begin{figure}[t]
	\begin{center}
		\figname{HMM_sketch}
		\newcommand{\dtLength}{4}
\newcommand{\arrowFactor}{0.33}
\newcommand{\slopes}{{0.3,0.05,-0.2}}

\newcommand{\mygreen}{green!60!black}
\begin{tikzpicture}[node distance = 2cm, auto]
	\tikzstyle{FHMMbarcur} = [-latex',line width=1,blue]
	\tikzstyle{FHMM1} = [-latex',line width=1,black]
	\tikzstyle{FHMM2} = [-latex',line width=1,\mygreen]
	\tikzstyle{FHMMnoise} =  [latex'-latex',dashed,line width=1,color=red]

	\draw [axis] (0,0) -- (13,0) node[right] {time};
	\draw [axis] (0,0) -- (0,3.5);
	
	\coordinate (P0) at (0.5,0.5);
	\xdef\yp {0.5}
	\foreach \iTime in {1, ..., 3}
	{
		\pgfmathparse{\yp + \slopes[\iTime-1]*\dtLength}
		\xdef\yp{\pgfmathresult}
		\coordinate (P\iTime) at ({\iTime*\dtLength+0.5},{\yp});
		
		\pgfmathtruncatemacro{\iTimePrev}{\iTime-1}
		\draw[dashed,line width=0.5pt] (P\iTimePrev) -- (P\iTime);
		
		\draw[FHMMbarcur] (P\iTimePrev) -- ++(\arrowFactor*\dtLength,{\slopes[\iTime-1]*\arrowFactor*\dtLength});
	}
	\draw[FHMM1] (P0) -- ++(\arrowFactor*\dtLength,{1.3*\arrowFactor*\dtLength}) node[above,shift={(1,0)}] {$\FHMM{\Xbar^{n-1}}{M}{\omega_n}$};
	\draw[FHMMnoise] ($(P0) + (\arrowFactor*\dtLength+0.05,{\slopes[0]*\arrowFactor*\dtLength})$) -- ++(0,{(1.3-\slopes[0])*\arrowFactor*\dtLength});
	\node[below,blue,shift={(0.6,-0.2)}] at ($(P0) + (\arrowFactor*\dtLength+0.05,{\slopes[0]*\arrowFactor*\dtLength})$) {$\FHMMbar{\Xbar^{n-1}}{M}$};
	
	\draw[FHMM1] (P1) -- ++(\arrowFactor*\dtLength,{1.1*\arrowFactor*\dtLength}) node[right,shift={(0.1,0)}] {$\FHMM{\Xbar^{n}}{M}{\omega_n}$};
	\draw[FHMMnoise] ($(P1) + (\arrowFactor*\dtLength+0.05,{\slopes[1]*\arrowFactor*\dtLength})$) -- ++(0,{(1.1-\slopes[1])*\arrowFactor*\dtLength});
	\node[right,blue,shift={(0.05,0.35)}] at ($(P1) + (\arrowFactor*\dtLength+0.05,{\slopes[1]*\arrowFactor*\dtLength})$) {$\FHMMbar{\Xbar^{n}}{M}$};
	\draw[FHMM2] (P1) -- ++(\arrowFactor*\dtLength,{-0.6*\arrowFactor*\dtLength});
	\draw[FHMMnoise] ($(P1) + (\arrowFactor*\dtLength+0.05,{\slopes[1]*\arrowFactor*\dtLength})$) -- ++(0,{(-0.6-\slopes[1])*\arrowFactor*\dtLength});
	\node[below,\mygreen,shift={(1,-0.1)}] at ($(P1) + (\arrowFactor*\dtLength+0.05,{-0.6*\arrowFactor*\dtLength})$) {$\FHMM{\Xbar^{n}}{M}{\omega_{n+1}}$};
	
	\draw[FHMM2] (P2) -- ++(\arrowFactor*\dtLength,{-0.8*\arrowFactor*\dtLength});
	\draw[FHMMnoise] ($(P2) + (\arrowFactor*\dtLength+0.05,{\slopes[2]*\arrowFactor*\dtLength})$) -- ++(0,{(-0.8-\slopes[2])*\arrowFactor*\dtLength});
	\node[right,blue,shift={(-0.2,0.35)}] at ($(P2) + (\arrowFactor*\dtLength+0.05,{\slopes[2]*\arrowFactor*\dtLength})$) {$\FHMMbar{\Xbar^{n+1}}{M}$};
	\node[below,\mygreen,shift={(1,-0.1)}] at ($(P2) + (\arrowFactor*\dtLength+0.05,{-0.8*\arrowFactor*\dtLength})$) {$\FHMM{\Xbar^{n+1}}{M}{\omega_{n+1}}$};

	\foreach \iTime in {0, ..., 3}
	{
		\def\tp {\iTime*\dtLength+0.5}
		\draw[dotted] (P\iTime) -- ({\tp},0); \draw ({\tp},-0.1) -- ({\tp},0.1);
		\draw[red,fill=red] (P\iTime) circle (0.4ex);
	}
	\node[right,shift={(-0.2,-0.35)}] at (0.5,0) {$t^{n-1}$};
	\node[right,shift={(-0.2,-0.35)}] at (\dtLength+0.5,0) {$t^{n}$};
	\node[right,shift={(-0.2,-0.35)}] at (2*\dtLength+0.5,0) {$t^{n+1}$};
	\node[right,shift={(-0.2,-0.35)}] at (3*\dtLength+0.5,0) {$t^{n+2}$};
	
	\path
	([shift={(-5\pgflinewidth,-5\pgflinewidth)}]current bounding box.south west)
	([shift={( 2\pgflinewidth, 2\pgflinewidth)}]current bounding box.north east);
\end{tikzpicture}
	\end{center}
	\vspace{-0.4cm}\caption{\label{fig:HMM_CV_sketch} Sketch of the proposed variance reduction technique. Blue arrows depict variance-reduced HMM estimators, which at every time instant require the calculation of two correlated original HMM estimators at the same and previous time instant. These original estimators are shown by black and green arrows for estimation at time $t^n$ and $t^{n+1}$, respectively. The red dashed arrows indicate that, due to the correlation, the estimator noise -- the difference between the black (or green) and blue arrows -- is similar at consecutive time instants. }
\end{figure}

The proposed technique can also be viewed from the following perspective. The difference between $\FHMM{\Xbar^{n-1}}{M}{\omega_n}$ and $\FHMMbar{\Xbar^{n-1}}{M}$ between brackets in equation \eqref{eq:F_cv} has zero expectation and approximately corresponds to the noise on the estimator $\FHMM{\Xbar^n}{M}{\omega_n}$ since the same seed $\omega_n$ is used in the first term. The variance reduction method is illustrated in figure \ref{fig:HMM_CV_sketch}.

\paragraph{Initialization} To get started, the procedure requires a variance-reduced estimation $\FHMMbar{\Xbar^0}{M}$ in the first step. There are several options:
\begin{itemize} 
	\item \emph{Exact solution.} In some of our numerical experiments, we will choose $\FHMMbar{\Xbar^0}{M}$ to be the exact solution $F(\Xbar^0)$. Clearly, this is a choice that cannot be made in practical applications (since it is not necessary to use the HMM method when this is possible), so this will only be done to illustrate some properties of the numerical scheme, most notably when studying the bias in section~\ref{subsec:est_bias}.
	\item \emph{More accurate HMM estimator.} A second option is to use a more accurately estimated value, denoted by $\FHMM{\Xbar^0}{M^*}{\omega_0}$, with a number of samples $M^*\gg M$. 
	\item \emph{Average of HMM estimators.} As a third option, one could also use an average of $S$ HMM estimators with $M$ realizations,
	\begin{equation}\label{eq:averaged_HMM_estimator}
		\FHMM{\Xbar^0}{S,M}{\boldsymbol{\omega}_0} = \frac{1}{S}\sum_{s=1}^S \FHMM{\Xbar^0}{M}{\omega_{0,s}},
	\end{equation}
	where $\boldsymbol{\omega}_0 = \set{\omega_{0,s}}{s}$ represents the vector of initial seeds. Notice that, when choosing $S=M^*/M$, computing~\eqref{eq:averaged_HMM_estimator} has the same computational cost as $\FHMM{\Xbar^0}{M^*}{\omega_0}$.
\end{itemize}
	
\paragraph{Reinitialization} Because the invariant measure that is sampled by the Markov chain~\eqref{eq:EM_scheme} is parametrized by the slow variable $X$, which itself evolves as a function of (macroscopic) time, we expect the variance reduction to become less effective as time advances, see also the analysis in section~\ref{subsec:est_var}. To reduce this effect, we introduce an additional reinitialization step: after every $R$ macroscopic time steps, we do not compute the variance reduced estimate as in~\eqref{eq:F_cv}, but instead use the initialization procedure described above.

\begin{rem}[Coarse projective integration] 
	The above procedure can also be used for the coarse projective integration estimator~\eqref{eq:F_cpi}, provided that the generation of the ensemble of initial conditions $\left\{y^{n,m}\right\}_{m=1}^M$ is done using the Euler-Maruyama scheme~\eqref{eq:EM_scheme}, and one ensures that the same seed $\omega_n$ is used when generating these initial conditions as well as for the Brownian increments to compute the time-evolved states $\left\{y^{n,m}_{K\delta t}\right\}_{m=1}^M$ in the two estimations. 
\end{rem}

\section{Numerical properties} \label{sec:num_prop}
The convergence of the HMM method described in section~\ref{subsec:hmm} has been studied in detail in the literature \cite{E2005,Fatkullin:JournalOfComputationalPhysics:2004}, see also \cite{Givon2006,Papavasiliou2007} for related results. In general, any HMM estimator contains errors from different sources.  First, while one intends to exploit the ergodicity property \eqref{eq:ergodicity_property}, one can only simulate the fast dynamics over a finite time interval $[0,\tau]$ with $\tau=M\dt$, which leads to a finite sampling error. 
Second, since the exact solution of the fast equation is not known explicitly, a time discretization method is used to approximate the solution of this equation, which leads to a discretization error. Third, we introduce an error by replacing the (finite) time integral by a finite Riemann sum, in which the fast variable is evaluated at discrete time instants, which leads to a sampling error. 

In the present paper, we are not concerned with these errors. We only study the reduction of the variance that results from superimposing the variance reduction technique of section~\ref{subsec:hmm_cv} onto the HMM estimator (section~\ref{subsec:est_var}). Subsequently, we study the potential \emph{additional} bias of the variance-reduced estimator~\eqref{eq:F_cv} with respect to the standard HMM estimator~\eqref{eq:F_hmm} in section~\ref{subsec:est_bias}.

\subsection{Estimator variance} \label{subsec:est_var}

\subsubsection{General case}
First, we consider the general (nonlinear) case and study the statistical error, which is quantified by the variance of the estimator:
\begin{equation} \label{eq:F_cv_variance}
	\Var[\FHMMbar{\Xbar^N}{M}] = \Var\left[\FHMM{\Xbar^N}{M}{\omega_N} - \left(\FHMM{\Xbar^{N-1}}{M}{\omega_N} - \FHMMbar{\Xbar^{N-1}}{M}\right)\right].
\end{equation}
Since the same Brownian path is used twice, we expect  $\FHMM{\Xbar^N}{M}{\omega_N}$ and $\FHMM{\Xbar^{N-1}}{M}{\omega_N}$ to be strongly correlated and the variance in the estimator $\FHMMbar{\Xbar^N}{M}$ reduced.
To obtain an expression for $\Var[\FHMMbar{\Xbar^N}{M}]$, we first rewrite the equation \eqref{eq:F_cv}. Starting from \eqref{eq:F_cv} and using \eqref{eq:F_hmm}, we get:
\begin{equation} \label{eq:F_diff}
	\FHMM{\Xbar^N}{M}{\omega_N} - \FHMM{\Xbar^{N-1}}{M}{\omega_N} = \frac{1}{M}\sum_{m=1}^M \left(f(\Xbar^N,y^{N,m}_N) - f(\Xbar^{N-1},y^{N-1,m}_N) \right),
\end{equation}
in which $\set{y^{n,m}_{\tilde{n}}}{m}$ denotes the set of Markov chain generated samples at time $t^n$ using a random number generator with seed $\omega_{\tilde{n}}$.
A Taylor expansion around $(\Xbar^N,y^{N,m}_N)$ of the difference within the summation of equation \eqref{eq:F_diff} leads to:
\begin{align} \label{eq:f_diff_taylor}
	f(\Xbar^N,y^{N,m}_N) - f(\Xbar^{N-1},y^{N-1,m}_N) &\approx \partial_xf^m\cdot(\Xbar^N - \Xbar^{N-1}) + \partial_yf^m\cdot(y^{N,m}_N - y^{N-1,m}_N) \notag \\
	&= \partial_xf^m\cdot\Dt\FHMMbar{\Xbar^{N-1}}{M} + \partial_yf^m\cdot(y^{N,m}_N - y^{N-1,m}_N),
\end{align}
where we used that fact that $\Xbar^N$ is obtained using a forward Euler step starting from $\Xbar^{N-1}$ (see equation~\eqref{eq:averaged_eq_FE_cv}) and introduced the shorthand notation $\partial_xf^m$ and $\partial_yf^m$ to denote the partial derivative of the function $f$ with respect to $x$ and $y$, respectively, evaluated at $(\Xbar^N,y^{N,m}_N)$ (since, for given $N$ the argument for which the partial derivative is evaluates is completely determined by $m$).
Next, the difference between the samples of the fast equation in the second term of equation \eqref{eq:f_diff_taylor} can be obtained by subtracting the Markov chains in equation \eqref{eq:EM_scheme} that generate them. For $m=1,\ldots,M-1$, this becomes:
\begin{align}
	y^{N,m+1}_N - y^{N-1,m+1}_N = y^{N,m}_N - y^{N-1,m}_N &+ \frac{\dt}{\epsi}(g(\Xbar^N,y^{N,m}_N) - g(\Xbar^{N-1},y^{N-1,m}_N)) \notag \\
	&+ \sqrt{\frac{\dt}{\epsi}}(\beta(\Xbar^N,y^{N,m}_N) - \beta(\Xbar^{N-1},y^{N-1,m}_N))\xi_N^m. \label{eq:y_diff}
\end{align}
We again use a Taylor expansion of the functions $g(x,y)$ and $\beta(x,y)$ around $(\Xbar^N,y^{N,m}_N)$ in \eqref{eq:y_diff}, yielding:
\begin{align}
	y^{N,m+1}_N - y^{N-1,m+1}_N \approx y^{N,m}_N - y^{N-1,m}_N &+ \frac{\dt}{\epsi}\left(\partial_x g^m\cdot(\Xbar^N - \Xbar^{N-1}) + \partial_y g^m\cdot(y^{N,m}_N - y^{N-1,m}_N)\right) \notag \\
	&+ \sqrt{\frac{\dt}{\epsi}}\left(\partial_x\beta^m\cdot(\Xbar^N - \Xbar^{N-1}) + \partial_y\beta^m\cdot(y^{N,m}_N - y^{N-1,m}_N)\right)\xi_N^m \label{eq:y_diff_taylor}.
\end{align}
Equation \eqref{eq:y_diff_taylor} can be compactly rewritten as:
\begin{equation} \label{eq:yvals_recursion}
	y^{N,m+1}_N - y^{N-1,m+1}_N \approx A_N^m(y^{N,m}_N - y^{N-1,m}_N) + B_N^m\Dt\FHMMbar{\Xbar^{N-1}}{M}, \quad m=1,...,M-1,
\end{equation}
in which the random numbers $A_N^m$ and $B_N^m$ are given by:
\begin{equation}\label{eq:coeffs_AB}
	A_N^m = 1 + \frac{\dt}{\epsi}\partial_yg^m + \sqrt{\frac{\dt}{\epsi}}\partial_y\beta^m\cdot\xi_N^m, \qquad\qquad 
	B_N^m = \frac{\dt}{\epsi}\partial_xg^m + \sqrt{\frac{\dt}{\epsi}}\partial_x\beta^m\cdot\xi_N^m.
\end{equation}
Working out equation \eqref{eq:yvals_recursion} leads to:
\begin{align}
	y^{N,1}_N - y^{N-1,1}_N &= 0 \notag \\
	y^{N,2}_N - y^{N-1,2}_N &\approx B_N^1\Dt\FHMMbar{\Xbar^{N-1}}{M} \notag \\
	y^{N,3}_N - y^{N-1,3}_N &\approx A_N^2(y^{N,2}_N - y^{N-1,2}_N) + B_N^2\Dt\FHMMbar{\Xbar^{N-1}}{M} \notag \\
							&= A_N^2B_N^1\Dt\FHMMbar{\Xbar^{N-1}}{M} + B_N^2\Dt\FHMMbar{\Xbar^{N-1}}{M} \notag \\
							&= (A_N^2B_N^1 + B_N^2)\Dt\FHMMbar{\Xbar^{N-1}}{M} \notag \\
	 y^{N,m}_N - y^{N-1,m}_N &\approx\sum_{i=1}^{m-1} \left(B_N^i\prod_{j=i+1}^{m-1} A_N^j \right)\Dt\FHMMbar{\Xbar^{N-1}}{M}, \qquad m \ge 2, \label{eq:yvals_norecursion}
\end{align} 
where we used in the first equation that both Markov chains start from the same initial condition $y^{N,1}_N = y^{N-1,1}_N = y^{N-1,M}_{N-1}$.
Substituting equation \eqref{eq:yvals_norecursion} into equation \eqref{eq:f_diff_taylor} we find:
\begin{equation*}
	f(\Xbar^N,y^{N,m}_N) - f(\Xbar^{N-1},y^{N-1,m}_N) \approx \left(\partial_xf^m + \partial_yf^m \sum_{i=1}^{m-1} \left(B_N^i \prod_{j=i+1}^{m-1} A_N^j \right)\right)\Dt\FHMMbar{\Xbar^{N-1}}{M},
\end{equation*}
from which we obtain:
\begin{equation} \label{eq:F_diff_approx}
	\FHMM{\Xbar^N}{M}{\omega_N} - \FHMM{\Xbar^{N-1}}{M}{\omega_N} \approx \frac{\Dt}{M}\sum_{m=1}^M\left(\partial_xf^m + \partial_yf^m \sum_{i=1}^{m-1} \left(B_N^i \prod_{j=i+1}^{m-1} A_N^j \right)\right)\FHMMbar{\Xbar^{N-1}}{M}.
\end{equation}
Substituting~\eqref{eq:F_diff_approx} into equation \eqref{eq:F_cv}, we find an approximation for the variance 
of the estimator in \eqref{eq:F_cv_variance}:
\begin{equation} \label{eq:F_cv_var}
	\Var\Big[\FHMMbar{\Xbar^N}{M}\Big] \approx \Var\Bigg[\FHMMbar{\Xbar^{N-1}}{M} + \frac{\Dt}{M}\sum_{m=1}^M\left(\partial_xf^m + \partial_yf^m \sum_{i=1}^{m-1} \left(B_N^i \prod_{j=i+1}^{m-1} A_N^j \right)\right)\FHMMbar{\Xbar^{N-1}}{M}\Bigg].
\end{equation}

From equation \eqref{eq:F_cv_var}, we can draw a number of conclusions: 
\begin{enumerate}[(i)]
	\item the variance of the estimator at time $t^N$ grows only slightly with respect to the variance at time $t^{N-1}$; 
	\item even when the variance at time $t^{N-1}$ is zero (which happens when it is computed via a deterministic reinitialization), the variance at time $t^N$  will be nonzero, since the coefficients $A^j_N$ and $B^j_N$ depend on the Brownian paths, see equation~\eqref{eq:coeffs_AB};
	\item the variance of the estimator will be an increasing function of $N$, since a bit of variance is added on at every macroscopic time step.
\end{enumerate}

This last observation is the reason we introduced a reinitialization procedure in section~\ref{subsec:hmm_cv}.

\subsubsection{The linear case}
For the linear system \eqref{eq:lin_slow_fast}, we show that the variance vanishes exactly when the initial estimator $\FHMMbar{X^0}{M}$ is deterministic. In this case, we have:
\begin{align*}
	&\partial_x f^m = \lambda ,\qquad \partial_y f^m = p ,\qquad \partial_xg^m = q ,\qquad \partial_yg^m = -A ,\qquad \partial_x\beta^m = \partial_y\beta^m = 0.
\end{align*}
Due to the linearity of the system, equation \eqref{eq:yvals_norecursion} (as well as all following equations) become exact. Moreover, because $\partial_x\beta^m = \partial_y\beta^m = 0$, the quantities $A_N^m$ and $B_N^m$, defined in~\eqref{eq:coeffs_AB}, become deterministic:
\begin{equation}
	A^m_N = 1 - A\frac{\dt}{\epsi}, \qquad B^m_N = q\frac{\dt}{\epsi}.
\end{equation}
Thus, equation~\eqref{eq:yvals_norecursion} can be rewritten as:
\begin{equation*}
	y^{N,m}_N - y^{N-1,m}_N = q\frac{\dt}{\epsi}\left(\sum_{i=0}^{m-2} \left(1 - A\frac{\dt}{\epsi}\right)^i\right)\Dt\FHMMbar{\Xbar^{N-1}}{M}, \qquad m \ge 2.
\end{equation*}
Equation \eqref{eq:F_diff_approx} reads:
\begin{equation} \label{eq:lin_F_diff}
	\FHMM{\Xbar^N}{M}{\omega_N} - \FHMM{\Xbar^{N-1}}{M}{\omega_N} = \left(\lambda + \frac{\dt}{\epsi}\frac{pq}{M}\sum_{m=2}^M \left(\sum_{i=0}^{m-2} \left(1 - A\frac{\dt}{\epsi}\right)^i\right)\right)\Dt\FHMMbar{\Xbar^{N-1}}{M},
\end{equation}
which is deterministic as soon as $\FHMMbar{\Xbar^{N-1}}{M}$ is deterministic. 
Using the known result on sums of geometric sequences, the sums in equation \eqref{eq:lin_F_diff} can be further calculated as:
\begin{equation} \label{eq:lin_F_diff_sums}
	\sum_{m=2}^M \left(\sum_{i=0}^{m-2} (1 - \Atilde)^i \right) = \frac{M}{\Atilde}\left(1 - \Btilde\right),
\end{equation}
where we introduced the following two constants:
\begin{equation} \label{eq:lin_Atilde_Btilde}
	\Atilde = A\dfrac{\dt}{\epsi}, \qquad \Btilde = \frac{1 - (1 - \Atilde)^{M}}{M\Atilde}.
\end{equation}
Combining equations \eqref{eq:lin_F_diff} and \eqref{eq:lin_F_diff_sums} and substituting the result into equation \eqref{eq:F_cv_without_brackets}, we find the following expression for the variance-reduced estimator in the linear case:
\begin{align}
	\FHMMbar{\Xbar^N}{M} &= \left( 1 + \Dt\left(\lambda + \frac{pq}{A}(1 - \Btilde)\right)\right)\FHMMbar{\Xbar^{N-1}}{M} \label{eq:lin_F_cv} \\
	& = \left( 1 + \Dt\left(\lambda + \frac{pq}{A}(1 - \Btilde)\right)\right)^N\FHMMbar{\Xbar^{0}}{M}. \label{eq:lin_F_cv_norecursion}
\end{align}
For the linear system, the expression of the variance in equation \eqref{eq:F_cv_var} is then given by:
\begin{equation} \label{eq:lin_F_var}
	\Var\left[\FHMMbar{\Xbar^N}{M}\right] = \left(1 + \Dt\left(\lambda + \frac{pq}{A}(1 - \Btilde)\right)\right)^{2N} \Var\left[\FHMMbar{\Xbar^0}{M}\right].
\end{equation}
Equation \eqref{eq:lin_F_var} reveals that in the linear case the variance of the proposed estimator at time step $t^N$ depends only on the variance of the initial estimation. If a deterministic initialization is used in the first forward Euler step in equation \eqref{eq:averaged_eq_FE_cv}, the variance vanishes since $\Var[\FHMMbar{\Xbar^0}{M}]=0$. 

This (perhaps) surprising result can also be seen as follows. In the second forward Euler step of the macroscopic equation, combining equation \eqref{eq:averaged_eq_FE_cv} with \eqref{eq:F_cv}, the variance of the variance-reduced estimator is readily obtained as:
\begin{align} \label{eq:var_F_cv_calc}
	\Var[\FHMMbar{\Xbar^{1}}{M}] &= \Var\left[\FHMM{\Xbar^1}{M}{\omega_1} - \left(\FHMM{\Xbar^0}{M}{\omega_1} - F(\Xbar^{0})\right)\right] \notag \\
	&= \Var\left[\frac{1}{M}\sum_{m=1}^M \left(\lambda (\Xbar^1-\Xbar^0) + p(y_1^{1,m}-y_1^{0,m})\right)\right] \notag \\
	&= \frac{p^2}{M^2}\Var\left[\sum_{m=1}^M \left(y_1^{1,m}-y_1^{0,m}\right)\right],
\end{align}
in which we used that $\Xbar^0$, $\Xbar^1$ and $\FHMMbar{\Xbar^0}{M}$ are all deterministic quantities. To calculate the sum in equation \eqref{eq:var_F_cv_calc}, we subtract the equations of the two Markov chains given in equation \eqref{eq:EM_scheme} from each other resulting in:
\begin{equation} \label{eq:lin_markov_chain_diff}
	y_1^{1,m+1}-y_1^{0,m+1} = \left(1 - A\frac{\dt}{\epsi}\right)(y_1^{1,m}-y_1^{0,m}) + \frac{\dt}{\epsi} q(\Xbar^1-\Xbar^0) + \sqrt{\frac{\dt}{\epsi}}(\xi_1^m-\xi_1^m).
\end{equation}
Since we are using the same Brownian path for both Markov chains, the stochastic part of equation \eqref{eq:lin_markov_chain_diff} cancels out exactly. Therefore, this difference between Markov chain generated samples is completely deterministic and its variance is zero. This continues to hold for all following forward Euler steps.

\subsection{Estimator bias} \label{subsec:est_bias}

\subsubsection{General case} \label{subsubsec:est_bias_general}
Next, we examine the bias of the variance-reduced HMM estimator at time instant $t^N$, with respect to the original HMM estimator.
Working out the recursion in equation \eqref{eq:F_cv}, the estimator can be written as follows:
\begin{alignat}{3}
	\FHMMbar{\Xbar^{N}}{M} &= \FHMM{\Xbar^N}{M}{\omega_N} &&- \left(\FHMM{\Xbar^{N-1}}{M}{\omega_N} - \FHMMbar{\Xbar^{N-1}}{M}\right) \label{eq:F_cv_without_brackets} \\
	&= \FHMM{\Xbar^N}{M}{\omega_N} &&+ \left(\FHMM{\Xbar^{N-1}}{M}{\omega_{N-1}} - \FHMM{\Xbar^{N-1}}{M}{\omega_N}\right) \notag \\
	& &&- \left(\FHMM{\Xbar^{N-2}}{M}{\omega_{N-1}} - \FHMMbar{\Xbar^{N-2}}{M}\right) \notag \\
	&= \FHMM{\Xbar^N}{M}{\omega_N} &&+ \sum_{n=1}^{N-1}\left(\FHMM{\Xbar^n}{M}{\omega_n} - \FHMM{\Xbar^n}{M}{\omega_{n+1}}\right) \notag \\
	& &&- \left(\FHMM{\Xbar^0}{M}{\omega_1} - \FHMMbar{\Xbar^0}{M}\right). \label{eq:F_cv_norecursion}
\end{alignat}
Taking the expectation of both sides of equation \eqref{eq:F_cv_norecursion} over repeated experiments while keeping the sample size $M$ and time step $\dt$ fixed leads to:
\begin{equation} \label{eq:F_cv_expectation}
	\E[\FHMMbar{\Xbar^{N}}{M}] = \E[\FHMM{\Xbar^N}{M}{\cdot}] - \left(\E[\FHMM{\Xbar^0}{M}{\cdot}] - \E[\FHMMbar{\Xbar^0}{M}]\right).
\end{equation}

From equation \eqref{eq:F_cv_expectation}, we observe that the variance-reduced HMM estimator $\FHMMbar{\Xbar^{N}}{M}$ does not introduce an additional bias compared to the original estimator $\FHMM{\Xbar^N}{M}{\omega_N}$, provided that the initial variance-reduced estimator $\FHMMbar{\Xbar^0}{M}$ is unbiased with respect to $\FHMM{\Xbar^0}{M}{\omega_0}$. One way of ensuring this is to use a Metropolized version of the HMM estimators, such that each individual term in equation~\eqref{eq:F_cv_expectation} is unbiased with respect to the exact time derivate $F(\cdot)$. (Note, however, that applying the Metropolis correction to both $\FHMM{\Xbar^N}{M}{\omega_N}$ and $\FHMM{\Xbar^{N-1}}{M}{\omega_N}$ may result in different samples getting rejected and hence a reduced correlation between the two estimator.)
Alternatively, one could try to ensure that the initial variance-reduced estimator $\FHMMbar{\Xbar^0}{M}$ contains exactly the same bias as $\FHMM{\Xbar^0}{M}{\omega_0}$. These effects are illustrated numerically in section~\ref{subsec:num_single_estimation}.

{\subsubsection{The linear case}
As equation~\eqref{eq:F_cv_expectation} shows, an additional bias may appear if the expectation of the initial estimator $\FHMMbar{\Xbar^0}{M}$ and of the standard HMM estimator $\FHMM{\Xbar^0}{M}{\omega_0}$ are different. We now calculate the resulting bias on the solution paths obtained with the variance-reduced HMM technique for the linear system~\eqref{eq:lin_slow_fast}. 
We first write the forward Euler solution of the averaged equation \eqref{eq:lin_averaged_eq} using the exact expression of $F$ given in equation \eqref{eq:lin_F_exact} as:
\begin{align} \label{eq:lin_X_fe_norecursion}
	X^{N} &= X^{N-1} + \Dt F(X^{N-1}) \notag
	= \left(1 + \Dt\left(\lambda + \frac{pq}{A}\right)\right)X^{N-1} \notag \\
	&= \left(1 + \Dt\left(\lambda + \frac{pq}{A}\right)\right)^N X^{0},
\end{align}
where the last line is obtained by working out the recursion.

We want to write a similar expression for the variance-reduced solution paths. To that end, using equation \eqref{eq:lin_F_cv_norecursion}, we write:
\begin{align} \label{eq:lin_X_cv_norecursion}
	\Xbar^{N} &= X^0 + \Dt \sum_{n=0}^{N-1}\FHMMbar{\Xbar^{n}}{M} \notag \\
	&= X^0 + \Dt\FHMMbar{\Xbar^{0}}{M} \sum_{n=0}^{N-1}\left( 1 + \Dt\left(\lambda + \frac{pq}{A}(1 - \Btilde)\right)\right)^n,
\end{align}
in which $\Btilde$ is defined in~\eqref{eq:lin_Atilde_Btilde} and the last line is again obtained by working out the recursion.
The sum in equation \eqref{eq:lin_X_cv_norecursion} corresponds to a geometric sequence and can be calculated as:
\begin{equation} \label{eq:lin_geometric_sum}
	\sum_{n=0}^{N-1}\left(1 + \Dt\left(\lambda + \frac{pq}{A}(1 - \Btilde)\right)\right)^n = \frac{\left(1 + \Dt\left(\lambda + \dfrac{pq}{A}(1 - \Btilde)\right)\right)^{N} - 1}{\Dt\left(\lambda + \dfrac{pq}{A}(1 - \Btilde)\right)}.
\end{equation}
Substituting this into equation \eqref{eq:lin_X_cv_norecursion} yields:
\begin{align} \label{eq:lin_X_cv_norecursion_part2}
	\Xbar^{N} = X^0 + \frac{\left(1 + \Dt\left(\lambda + \dfrac{pq}{A}(1 - \Btilde)\right)\right)^{N} - 1}{\left(\lambda + \dfrac{pq}{A}(1 - \Btilde)\right)} \FHMMbar{\Xbar^{0}}{M}.
\end{align}
Equation \eqref{eq:lin_X_cv_norecursion_part2} provides the variance-reduced solution path at time $t^N$ which depends solely on the initial estimation $\FHMMbar{\Xbar^{0}}{M}$. It is straightforward to verify that this formula delivers the expected results for $N=0$ and $N=1$. In addition, equation \eqref{eq:lin_X_cv_norecursion_part2} allows to calculate the asymptotic $(N \to \infty)$ behavior of $\Xbar^N$ for fixed $\Dt$. 
Since $M > 0$ and $\Atilde = A\dfrac{\dt}{\epsi} \le 2$ due to stability of the Euler-Maruyama scheme in equation \eqref{eq:EM_scheme}, 
we obtain:
\begin{equation}
 	1 - \Btilde \le 1 < - \frac{\lambda A}{pq},
 \end{equation}
 with $\Btilde$ defined in~\eqref{eq:lin_Atilde_Btilde}, and hence equation \eqref{eq:lin_X_cv_norecursion_part2} yields:
\begin{equation} \label{eq:lin_X_cv_bias}
	\lim_{N\to\infty} \Xbar^N = \Xbar^{\infty} = \Xbar^0 - \frac{\FHMMbar{\Xbar^{0}}{M}}{\lambda + \dfrac{pq}{A}(1 - \Btilde)}.
\end{equation}
When using an exact initialization $\FHMMbar{\Xbar^{0}}{M} = F(X^0)$, we find:
\begin{align} \label{eq:lin_X_cv_bias_exact}
	\Xbar^{\infty} = \Xbar^0 - \frac{\left(\lambda + \dfrac{pq}{A}\right)\Xbar^{0}}{\lambda + \dfrac{pq}{A}(1 - \Btilde)} = - \frac{\dfrac{pq}{A}\Btilde\Xbar^{0}}{\left(\lambda + \dfrac{pq}{A}(1 - \Btilde)\right)}
\end{align}
From equation \eqref{eq:lin_X_cv_bias_exact} we learn that the factor $\Btilde$ determines the additional bias in solution paths of the linear system. To avoid a bias compared to the forward Euler solution with exact $F$ we require that $\Btilde = 0$. Using the expression of $\Btilde$ in equation \eqref{eq:lin_Atilde_Btilde} the following condition arises:
\begin{equation*}
	\Btilde = 0 \Rightarrow \left(1 - A\dfrac{\dt}{\epsi}\right)^{M} = 1,
\end{equation*}
which is satisfied in the following three cases: (i) $M=0$, (ii) $\epsi$ constant and $\dt \to 0$, and (iii) if the sample size $M$ is even and at the same time the following constraint holds:
\begin{equation} \label{eq:lin_nobias_condition}
	A\dfrac{\dt}{\epsi} = 2.
\end{equation}
In the numerical experiments in section~\ref{sec:results}, we only consider the third case.

\section{Numerical results} \label{sec:results}
We now put the variance reduction method to the test. We begin by looking at the bias and variance of the method after one iteration in section \ref{subsec:num_single_estimation}. We then turn to more detailed experiments on local variance reduction as a function of the numerical parameters in section~\ref{subsec:local_var_reduction} and on the resulting variance on the solution trajectories in \ref{subsec:trajectories}.  In all cases, we compare the linear and nonlinear model problems~\eqref{eq:lin_slow_fast} and~\eqref{eq:nonlin_slow_fast}.

\subsection{Bias and variance of a single estimation} \label{subsec:num_single_estimation}
Here, we investigate the effect of the three possible initialization procedures described in section~\ref{subsec:hmm_cv}. Since the first forward Euler step in any variance-reduced scheme needs to be taken with an accurate (low-variance) initial estimate $\FHMMbar{X^0}{M}$, we compare all estimations at time $t^1 = \Dt$, given the macroscopic state $\Xbar^1$, obtained as:
\begin{equation} \label{eq:bias_var_X1}
	\Xbar^1 = X^0 + \Dt\FHMMbar{X^0}{M}.
\end{equation}
We perform $J_r$ realizations, denoted by $(\Xbar_j^1)_{j=1}^{J_r}$, using different random seeds.
In section \ref{subsubsec:bias_var_lin}, we discuss the linear model problem~\eqref{eq:lin_slow_fast}. Afterwards, in section~\ref{subsubsec:bias_var_nonlin}, we consider the nonlinear system given in \eqref{eq:nonlin_slow_fast}. In both sections, we set $\Dt = 0.02$, $J_r = 500$ and $\epsi = 10^{-3}$.

\subsubsection{Linear system} \label{subsubsec:bias_var_lin}
For the linear system \eqref{eq:lin_slow_fast}, we choose $\lambda = -10$, $p = 4$, $q = 0.5$ and $A = 1.2$. The initial conditions are chosen as $X^0 = x_0 = 1$ and $y_0 = 1$. To show the variance reduction of the proposed method, we also perform a reference simulation using HMM without variance reduction. For this simulation, we choose the Markov chain time step $\dt = \epsi$, while for the variance-reduced estimator, we fix $\dt$ as given in equation \eqref{eq:lin_nobias_condition}.

\paragraph{HMM without variance reduction} We begin by computing the HMM estimator distribution without variance reduction at time $t^1$ by calculating 500 realizations $\FHMM{\Xhat^1_j}{M}{\omega_{1,j}}$ using $M = 50$ samples and different seeds $(\omega_{1,j})_{j=1}^{J_r}$. Recall that the samples are generated by simulating the fast equation of system \eqref{eq:lin_slow_fast} using the Euler-Maruyama scheme given in equation \eqref{eq:EM_scheme}, while keeping the current value of the slow variable fixed. The first forward Euler step is performed similar to \eqref{eq:bias_var_X1}:
\begin{equation} \label{eq:bias_var_X1_HMM}
	\Xhat^1_j = X^0 + \Dt\FHMM{X^0}{M}{\omega_{0,j}}, \qquad 1 \le j \le J_r,
\end{equation}
with different seeds $(\omega_{0,j})_{j=1}^{J_r}$ in each realization.
We compare the difference between estimations with and without the Metropolis-Hastings correction.
In figure \ref{fig:pdfs_linear} (left), we plot the estimator distributions at time $t^1$, in which the blue and green dot-dashed lines represent the distributions without and with the Metropolis-Hastings algorithm, respectively. The vertical lines of the same color depict the sample mean of each distribution. The vertical red line corresponds to the exact mean, given by $F(X^1)$, with $X^1$ the forward Euler solution of the averaged equation using the exact expression of $F(x)$ in equation \eqref{eq:lin_F_exact}. We observe that, in the linear case, the HMM estimator is unbiased with respect to the forward Euler solution for $M \to \infty$. (This is due to the fact that the Markov chain generated by the Euler-Maruyama scheme~\eqref{eq:EM_scheme} preserves the mean of the fast variable $y$, which is the only information on $y$ that is used in the effective equation~\eqref{eq:lin_averaged_eq}.) Moreover, the addition of the extra Metropolis step has no significant influence on the variance of the estimator $\FHMM{\Xbar^1}{M}{\omega_1}$.

\begin{figure}[t]
	\begin{center}
		\figname{pdfs_linear}
%
%
\definecolor{mycolor1}{rgb}{0.00000,0.49804,0.00000}%
\definecolor{mycolor2}{rgb}{1.00000,0.00000,1.00000}%
\newcommand{\figWidth}{0.27\textwidth} 
\newcommand{\figHeight}{3cm} 
\newcommand{\figSpacingRight}{1cm} 
\begin{tikzpicture}

\begin{axis}[%
width=0,
height=0,
scale only axis,
clip=false,
xmin=0,
xmax=1,
ymin=0,
ymax=1,
hide axis,
]
\node[align=center, text=black]
at (3*\figWidth/2+\figSpacingRight,\figHeight+0.9cm) {Estimator distributions (linear case)};
\end{axis}

\begin{axis}[%
width=\figWidth,
height=\figHeight,
at={(0,0.0)},
scale only axis,
xmin=-8.5,
xmax=-5,
ymin=0,
ymax=3.5,
axis background/.style={fill=white},
title style={font=\scriptsize\bfseries},
title={no variance reduction},
axis x line*=bottom,
axis y line*=left,
legend style={shift={(0,0)},legend cell align=left,align=left,draw=white!15!black}
]
\addplot [color=blue,dashdotted]
  table[]{tikz/data/linear_pdfs-1.tsv};
\addlegendentry{HMM};

\addplot [color=mycolor1,dashdotted]
  table[]{tikz/data/linear_pdfs-2.tsv};
\addlegendentry{HMM + MH};

\addplot [color=red,solid]
  table[]{tikz/data/linear_pdfs-3.tsv};
\addlegendentry{Exact mean};

\addplot [color=blue,solid,forget plot]
  table[]{tikz/data/linear_pdfs-4.tsv};
\addplot [color=mycolor1,solid,forget plot]
  table[]{tikz/data/linear_pdfs-5.tsv};
\end{axis}

\begin{axis}[%
width=\figWidth,
height=\figHeight,
at={(\figWidth+\figSpacingRight,0)},
scale only axis,
xmin=-8.5,
xmax=-5,
ymin=0,
ymax=3.5,
axis background/.style={fill=white},
title style={font=\scriptsize\bfseries},
title={variance reduction without MH},
axis x line*=bottom,
axis y line*=left,
legend style={shift={(0.1,0)},legend cell align=left,align=left,draw=white!15!black}
]
\addplot [color=mycolor2,solid]
  table[]{tikz/data/linear_pdfs-6.tsv};
\addlegendentry{Exact};

\addplot [color=blue,dashdotted]
  table[]{tikz/data/linear_pdfs-7.tsv};
\addlegendentry{Estimated};

\addplot [color=mycolor1,dashdotted]
  table[]{tikz/data/linear_pdfs-8.tsv};
\addlegendentry{Averaged};

\addplot [color=red,solid]
  table[]{tikz/data/linear_pdfs-9.tsv};
\addlegendentry{Exact mean};

\addplot [color=blue,solid,forget plot]
  table[]{tikz/data/linear_pdfs-10.tsv};
\addplot [color=mycolor1,solid,forget plot]
  table[]{tikz/data/linear_pdfs-11.tsv};
\end{axis}

\begin{axis}[%
width=\figWidth,
height=\figHeight,
at={(2*\figWidth+2*\figSpacingRight,0)},
scale only axis,
xmin=-8.5,
xmax=-5,
ymin=0,
ymax=3.5,
axis background/.style={fill=white},
title style={font=\scriptsize\bfseries},
title={variance reduction with MH},
axis x line*=bottom,
axis y line*=left,
legend style={shift={(0.1,0)},legend cell align=left,align=left,draw=white!15!black}
]
\addplot [color=blue,dashdotted]
  table[]{tikz/data/linear_pdfs-12.tsv};
\addlegendentry{Estimated};

\addplot [color=mycolor1,dashdotted]
  table[]{tikz/data/linear_pdfs-13.tsv};
\addlegendentry{Averaged};

\addplot [color=red,solid]
  table[]{tikz/data/linear_pdfs-14.tsv};
\addlegendentry{Exact mean};

\addplot [color=blue,solid,forget plot]
  table[]{tikz/data/linear_pdfs-15.tsv};
\addplot [color=mycolor1,solid,forget plot]
  table[]{tikz/data/linear_pdfs-16.tsv};
\end{axis}
\end{tikzpicture}%
	\end{center}
	\caption{\label{fig:pdfs_linear} Estimator distributions after one time step $\Dt = 0.02$ when applied to the linear system \eqref{eq:lin_slow_fast}. The exact mean $F(\Xbar^1)$ is depicted by a vertical red line in each plot. Left: HMM without variance reduction with (green) and without (blue) Metropolis-Hastings. Middle: variance-reduced HMM without Metropolis-Hastings for different initializations. Right: variance-reduced HMM with Metropolis-Hastings used in the estimated and averaged initializations. }
\end{figure}

\paragraph{HMM with variance reduction without Metropolis-Hastings} Next, we calculate the variance-reduced estimator distribution at time $t^1$ and examine the effect of the different possible initializations (choices for $\FHMMbar{X^0}{M}$). In the middle plot of figure \ref{fig:pdfs_linear}, we show the results without adding the Metropolis-Hastings correction at $t^0$. The blue and green dot-dashed lines correspond to the distributions when using an estimated initialization $\FHMM{X^0}{500}{\omega_{0,j}}$ and an averaged initialization $\FHMM{X^0}{10,50}{\omega_{0,j}}$ in each realization $j=1,...,J_r$, respectively. The pink line represents the estimator distribution when using the exact expression $F(X^0)$ as initial estimator. Since the variance-reduced estimator yields a variance-free result when choosing an exact initialization in the linear case, see equation \eqref{eq:lin_F_var}, the corresponding distribution reduces to a vertical pink line. The vertical red line corresponds to the exact mean $F(X^1)$. Since all mean values lie close together, the pink and red lines are hard to discern. We find that both initializations lead to an unbiased variance-reduced estimator and the reduction in variance is clearly visible.
In figure \ref{fig:pdfs_linear_est_av}, we more closely inspect the estimated and averaged initializations. On the left, we show the estimator distributions when using $\FHMM{X^0}{500}{\omega_{0,j}}$ (solid blue) and $\FHMM{X^0}{5000}{\omega_{0,j}}$ (dashed blue) as initial estimators. This confirms that the estimated initialization indeed leads to an unbiased estimator for $M^* \to \infty$ and yields reductions in variance by a factor $12$ and $118$, respectively. 
On the right, we visualize the distributions when using $\FHMM{X^0}{10,50}{\omega_{0,j}}$ (solid green), $\FHMM{X^0}{100,50}{\omega_{0,j}}$ (dashed green) and $\FHMM{X^0}{1000,50}{\omega_{0,j}}$ (dot-dashed green). This shows that the averaged initialization also gives rise to an unbiased estimator for $S \to \infty$ with $M$ fixed and yields reductions by a factor $13$, $115$ and $1163$, respectively.

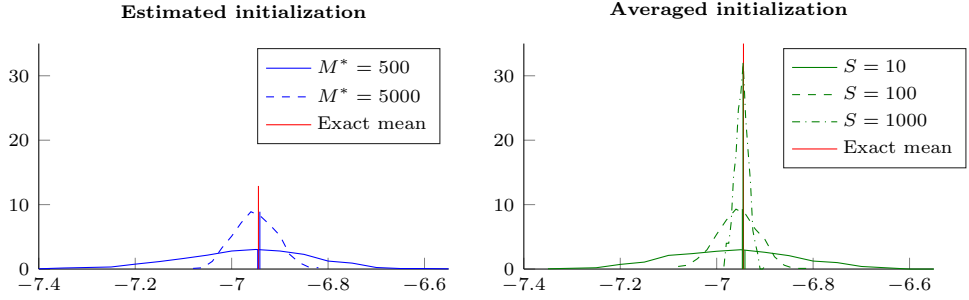
\begin{figure}[t]
	\begin{center}
		\figname{pdfs_linear_est_av}
%
%
\definecolor{mycolor1}{rgb}{0.00000,0.49804,0.00000}%
\newcommand{\figWidth}{0.33\textwidth} 
\newcommand{\figHeight}{3cm} 
\newcommand{\figSpacingRight}{1cm} 
\begin{tikzpicture}

\begin{axis}[%
width=0,
height=0,
scale only axis,
clip=false,
xmin=0,
xmax=1,
ymin=0,
ymax=1,
hide axis,
]
\node[align=center, text=black]
at (\figWidth+\figSpacingRight/2,\figHeight+0.9cm) {Variance-reduced estimator distributions without MH (linear case)};
\end{axis}

\begin{axis}[%
width=\figWidth,
height=\figHeight,
at={(0,0.0)},
scale only axis,
xmin=-7.4,
xmax=-6.55,
ymin=0,
ymax=35,
axis background/.style={fill=white},
title style={font=\scriptsize\bfseries},
title={Estimated initialization},
axis x line*=bottom,
axis y line*=left,
legend style={legend cell align=left,align=left,draw=white!15!black}
]
\addplot [color=blue,solid]
  table[]{tikz/data/linear_est_av_pdfs-1.tsv};
\addlegendentry{$M^{*} = 500$};

\addplot [color=blue,dashed]
  table[]{tikz/data/linear_est_av_pdfs-2.tsv};
\addlegendentry{$M^{*} = 5000$};

\addplot [color=red,solid]
  table[]{tikz/data/linear_est_av_pdfs-4.tsv};
\addlegendentry{Exact mean};

\addplot [color=blue,solid,forget plot]
  table[]{tikz/data/linear_est_av_pdfs-5.tsv};
\addplot [color=blue,solid,forget plot]
  table[]{tikz/data/linear_est_av_pdfs-6.tsv};
\end{axis}

\begin{axis}[%
width=\figWidth,
height=\figHeight,
at={(\figWidth+\figSpacingRight,0)},
scale only axis,
xmin=-7.4,
xmax=-6.55,
ymin=0,
ymax=35,
axis background/.style={fill=white},
title style={font=\scriptsize\bfseries},
title={Averaged initialization},
axis x line*=bottom,
axis y line*=left,
legend style={shift={(0.1,0)},legend cell align=left,align=left,draw=white!15!black}
]
\addplot [color=mycolor1,solid]
  table[]{tikz/data/linear_est_av_pdfs-7.tsv};
\addlegendentry{$S = 10$};

\addplot [color=mycolor1,dashed]
  table[]{tikz/data/linear_est_av_pdfs-8.tsv};
\addlegendentry{$S = 100$};

\addplot [color=mycolor1,dashdotted]
  table[]{tikz/data/linear_est_av_pdfs-9.tsv};
\addlegendentry{$S = 1000$};

\addplot [color=red,solid]
  table[]{tikz/data/linear_est_av_pdfs-11.tsv};
\addlegendentry{Exact mean};

\addplot [color=mycolor1,solid,forget plot]
  table[]{tikz/data/linear_est_av_pdfs-12.tsv};
\addplot [color=mycolor1,solid,forget plot]
  table[]{tikz/data/linear_est_av_pdfs-13.tsv};
\addplot [color=mycolor1,solid,forget plot]
  table[]{tikz/data/linear_est_av_pdfs-14.tsv};
\end{axis}
\end{tikzpicture}%
	\end{center}
	\caption{\label{fig:pdfs_linear_est_av} Variance-reduced estimator distributions (without Metropolis-Hastings) after one time step $\Dt = 0.02$ when applied to the linear system \eqref{eq:lin_slow_fast}. The exact mean $F(X^1)$ is depicted by a vertical red line in each plot. Left: estimated initialization. Right: averaged initialization. }
\end{figure}

\paragraph{HMM with variance reduction with Metropolis-Hastings} When using Metropolis-Hastings to generate the samples in the initialization $\FHMMbar{X^0}{M}$, we obtain the distributions in figure \ref{fig:pdfs_linear} (right). The effect of using $\FHMM{X^0}{500}{\omega_{0,j}}$ and $\FHMM{X^0}{10,50}{\omega_{0,j}}$ as initial estimators is shown by blue and green dot-dashed lines, respectively. The former leads to an unbiased estimator with clear variance reduction, while the latter results in a biased estimator with only little reduction in variance. 
We regard both initializations in more detail in figure \ref{fig:pdfs_linear_est_av_MH}. By comparing the effect of using $\FHMM{X^0}{500}{\omega_{0,j}}$ (blue solid) and $\FHMM{X^0}{5000}{\omega_{0,j}}$ (blue dashed) in the left plot, we conclude that the estimated initialization yields an unbiased estimator for $M^* \to \infty$ with reduction factors factors $11$ and $106$. The right plot indicates that by using Metropolis-Hastings in the averaged initialization the resulting estimator becomes unbiased for $S \to \infty$ with $M$ fixed. Moreover, the reduction in variance is significantly lower with corresponding reduction factors $1$, $2$ and $21$. Although both initializations $\FHMM{X^0}{M^*}{\omega_{0,j}}$ and $\FHMM{X^0}{M,S}{\omega_{0,j}}$ are computationally equivalent when $M^* = SM$, the former outperforms the latter by far, which requires further research to better understand this behavior.

\begin{figure}[t]
	\begin{center}
		\figname{pdfs_linear_est_av_MH}
%
%
\definecolor{mycolor1}{rgb}{0.00000,0.49804,0.00000}%
\newcommand{\figWidth}{0.33\textwidth} 
\newcommand{\figHeight}{3cm} 
\newcommand{\figSpacingRight}{1cm} 
\begin{tikzpicture}

\begin{axis}[%
width=0,
height=0,
scale only axis,
clip=false,
xmin=0,
xmax=1,
ymin=0,
ymax=1,
hide axis,
]
\node[align=center, text=black]
at (\figWidth+\figSpacingRight/2,\figHeight+0.9cm) {Variance-reduced estimator distributions with MH (linear case)};
\end{axis}

\begin{axis}[%
width=\figWidth,
height=\figHeight,
at={(0,0.0)},
scale only axis,
xmin=-8,
xmax=-5.57,
ymin=0,
ymax=11,
ytick={0, 2, ..., 14},
axis background/.style={fill=white},
title style={font=\scriptsize\bfseries},
title={Estimated initialization},
axis x line*=bottom,
axis y line*=left,
legend style={legend cell align=left,align=left,draw=white!15!black}
]
\addplot [color=blue,solid]
  table[]{tikz/data/linear_est_av_MH_pdfs-1.tsv};
\addlegendentry{$M^{*} = 500$};

\addplot [color=blue,dashed]
  table[]{tikz/data/linear_est_av_MH_pdfs-2.tsv};
\addlegendentry{$M^{*} = 5000$};

\addplot [color=red,solid]
  table[]{tikz/data/linear_est_av_MH_pdfs-5.tsv};
\addlegendentry{Exact mean};

\addplot [color=blue,solid,forget plot]
  table[]{tikz/data/linear_est_av_MH_pdfs-3.tsv};
\addplot [color=blue,solid,forget plot]
  table[]{tikz/data/linear_est_av_MH_pdfs-4.tsv};
\end{axis}

\begin{axis}[%
width=\figWidth,
height=\figHeight,
at={(\figWidth+\figSpacingRight,0)},
scale only axis,
xmin=-8,
xmax=-5.57,
ymin=0,
ymax=11,
ytick={0, 2, ..., 14},
axis background/.style={fill=white},
title style={font=\scriptsize\bfseries},
title={Averaged initialization},
axis x line*=bottom,
axis y line*=left,
legend style={legend cell align=left,align=left,draw=white!15!black}
]
\addplot [color=mycolor1,solid]
  table[]{tikz/data/linear_est_av_MH_pdfs-6.tsv};
\addlegendentry{$S = 10$};

\addplot [color=mycolor1,dashed]
  table[]{tikz/data/linear_est_av_MH_pdfs-7.tsv};
\addlegendentry{$S = 100$};

\addplot [color=mycolor1,dashdotted]
  table[]{tikz/data/linear_est_av_MH_pdfs-8.tsv};
\addlegendentry{$S = 1000$};

\addplot [color=red,solid]
  table[]{tikz/data/linear_est_av_MH_pdfs-12.tsv};
\addlegendentry{Exact mean};

\addplot [color=mycolor1,solid,forget plot]
  table[]{tikz/data/linear_est_av_MH_pdfs-9.tsv};
\addplot [color=mycolor1,solid,forget plot]
  table[]{tikz/data/linear_est_av_MH_pdfs-10.tsv};
\addplot [color=mycolor1,solid,forget plot]
  table[]{tikz/data/linear_est_av_MH_pdfs-11.tsv};
\end{axis}
\end{tikzpicture}%
	\end{center}
	\vspace{-0.4cm}\caption{\label{fig:pdfs_linear_est_av_MH} Variance-reduced estimator distributions (including Metropolis-Hastings) after one time step $\Dt = 0.02$ when applied to the linear system \eqref{eq:lin_slow_fast}. The exact mean $F(X^1)$ is depicted by a vertical red line in each plot. Left: estimated initialization. Right: averaged initialization. }
\end{figure}
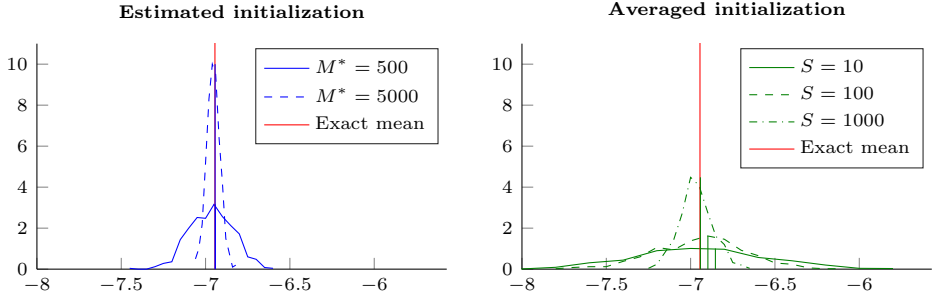

\subsubsection{Nonlinear system} \label{subsubsec:bias_var_nonlin}
We now consider the nonlinear system \eqref{eq:nonlin_slow_fast}, and set the initial conditions as $x_0 = 1$ and $y_0 = 0.5$. Furthermore, we set the time step $\dt = \epsi$. In the linear case, we compared our results with the exact mean at time $t^1$, given by $F(X^1)$. However, since the HMM estimator is biased in the general nonlinear case and we are only interested in studying the bias in the estimation (and not in the solution paths $\Xbar^1$), here, we define the exact mean for the HMM estimator without and with variance reduction as:
\begin{equation} \label{eq:bias_var_nonlin_mean_exact}
	\widehat{m}_e = \frac{1}{J_r}\sum_{j=1}^{J_r} F(\Xhat^1_j), \qquad \bar{m}_e = \frac{1}{J_r}\sum_{j=1}^{J_r} F(\Xbar^1_j),
\end{equation}
in which $F(x)$ corresponds to the exact right hand side of the averaged equation given in equation \eqref{eq:nonlin_F_exact}.

\paragraph{HMM without variance reduction} When computing the HMM estimator distribution at time $t^1$ without and with Metropolis-Hastings, we obtain the blue and green dot-dashed distributions in figure \ref{fig:pdfs_nonlinear} (left), respectively. The vertical red lines correspond to the exact means $\widehat{m}_e$ given in equation \eqref{eq:bias_var_nonlin_mean_exact} for both experiments, which coincide to the naked eye. We immediately see the necessity of including the Metropolis-Hastings algorithm: the blue distribution has a very clear bias after only one iteration of the method, while the green distribution possesses the correct mean.

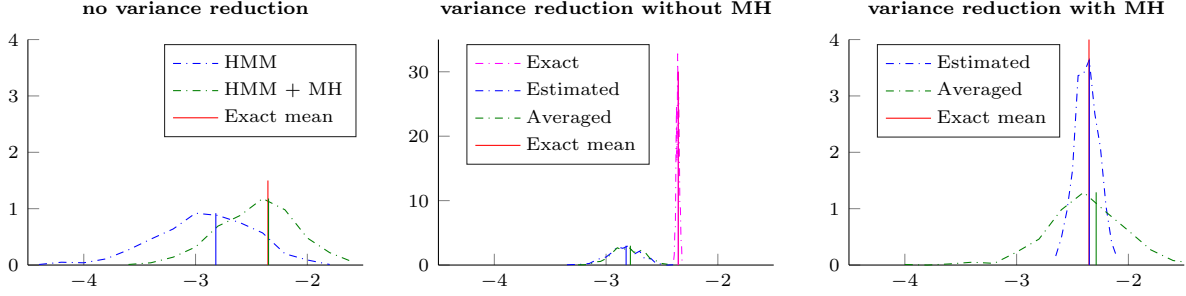
\begin{figure}[t]
	\begin{center}
		\figname{pdfs_nonlinear}
%
%
\definecolor{mycolor1}{rgb}{0.00000,0.49804,0.00000}%
\definecolor{mycolor2}{rgb}{1.00000,0.00000,1.00000}%
\newcommand{\figWidth}{0.27\textwidth} 
\newcommand{\figHeight}{3cm} 
\newcommand{\figSpacingRight}{1cm} 
\begin{tikzpicture}

\begin{axis}[%
width=0,
height=0,
scale only axis,
clip=false,
xmin=0,
xmax=1,
ymin=0,
ymax=1,
hide axis,
]
\node[align=center, text=black]
at (3*\figWidth/2+\figSpacingRight,\figHeight+0.9cm) {Estimator distributions (nonlinear case)};
\end{axis}

\begin{axis}[%
width=\figWidth,
height=\figHeight,
at={(0,0.0)},
scale only axis,
xmin=-4.5,
xmax=-1.5,
ymin=0,
ymax=4,
axis background/.style={fill=white},
title style={font=\scriptsize\bfseries},
title={no variance reduction},
axis x line*=bottom,
axis y line*=left,
legend style={legend cell align=left,align=left,draw=white!15!black}
]
\addplot [color=blue,dashdotted]
  table[]{tikz/data/nonlinear_pdfs-1.tsv};
\addlegendentry{HMM};

\addplot [color=mycolor1,dashdotted]
  table[]{tikz/data/nonlinear_pdfs-2.tsv};
\addlegendentry{HMM + MH};

\addplot [color=red,solid]
  table[]{tikz/data/nonlinear_pdfs-5.tsv};
\addlegendentry{Exact mean};

\addplot [color=mycolor1,solid,forget plot]
  table[]{tikz/data/nonlinear_pdfs-3.tsv};
\addplot [color=blue,solid,forget plot]
  table[]{tikz/data/nonlinear_pdfs-4.tsv};
\end{axis}

\begin{axis}[%
width=\figWidth,
height=\figHeight,
at={(\figWidth+\figSpacingRight,0)},
scale only axis,
xmin=-4.5,
xmax=-1.5,
ymin=0,
ymax=35,
axis background/.style={fill=white},
title style={font=\scriptsize\bfseries},
title={variance reduction without MH},
axis x line*=bottom,
axis y line*=left,
legend style={shift={(-0.35,0)},legend cell align=left,align=left,draw=white!15!black}
]
\addplot [color=mycolor2,dashdotted]
  table[]{tikz/data/nonlinear_pdfs-6.tsv};
\addlegendentry{Exact};

\addplot [color=blue,dashdotted]
  table[]{tikz/data/nonlinear_pdfs-7.tsv};
\addlegendentry{Estimated};

\addplot [color=mycolor1,dashdotted]
  table[]{tikz/data/nonlinear_pdfs-8.tsv};
\addlegendentry{Averaged};

\addplot [color=red,solid]
  table[]{tikz/data/nonlinear_pdfs-12.tsv};
\addlegendentry{Exact mean};

\addplot [color=mycolor2,solid,forget plot]
  table[]{tikz/data/nonlinear_pdfs-9.tsv};
\addplot [color=mycolor1,solid,forget plot]
  table[]{tikz/data/nonlinear_pdfs-10.tsv};
\addplot [color=blue,solid,forget plot]
  table[]{tikz/data/nonlinear_pdfs-11.tsv};
\end{axis}

\begin{axis}[%
width=\figWidth,
height=\figHeight,
at={(2*\figWidth+2*\figSpacingRight,0)},
scale only axis,
xmin=-4.5,
xmax=-1.5,
ymin=0,
ymax=4,
axis background/.style={fill=white},
title style={font=\scriptsize\bfseries},
title={variance reduction with MH},
axis x line*=bottom,
axis y line*=left,
legend style={shift={(-0.35,0)},legend cell align=left,align=left,draw=white!15!black}
]
\addplot [color=blue,dashdotted]
  table[]{tikz/data/nonlinear_pdfs-13.tsv};
\addlegendentry{Estimated};

\addplot [color=mycolor1,dashdotted]
  table[]{tikz/data/nonlinear_pdfs-14.tsv};
\addlegendentry{Averaged};

\addplot [color=red,solid]
  table[]{tikz/data/nonlinear_pdfs-16.tsv};
\addlegendentry{Exact mean};

\addplot [color=mycolor1,solid,forget plot]
  table[]{tikz/data/nonlinear_pdfs-15.tsv};
\addplot [color=blue,solid,forget plot]
  table[]{tikz/data/nonlinear_pdfs-17.tsv};
\end{axis}
\end{tikzpicture}%
	\end{center}
	\caption{\label{fig:pdfs_nonlinear} Estimator distributions after one time step $\Dt = 0.02$ when applied to the nonlinear system \eqref{eq:nonlin_slow_fast}. The exact mean given in \eqref{eq:bias_var_nonlin_mean_exact} is depicted by a vertical red line in each plot. Left: HMM without variance reduction with (green) and without (blue) Metropolis-Hastings. Middle: variance-reduced HMM without Metropolis-Hastings for different initializations. Right: variance-reduced HMM with Metropolis-Hastings used in the estimated and averaged initializations. }
\end{figure}

\paragraph{HMM with variance reduction without Metropolis-Hastings} Repeating the above experiment for the variance-reduced estimator without Metropolis-Hastings at time $t^1$ using different initializations, we show the results in the middle plot of figure \ref{fig:pdfs_nonlinear}. The blue and green dot-dashed lines represent the distributions when using $\FHMM{X^0}{500}{\omega_{0,j}}$ and an $\FHMM{X^0}{10,50}{\omega_{0,j}}$ in each realization $j=1,...,J_r$, respectively. When using an exact initialization $F(X^0)$, we obtain the sharply peaked pink distribution centered around the exact means $\bar{m}_e$ (vertical red lines) in equation \eqref{eq:bias_var_nonlin_mean_exact} of these three experiments, which again coincide to the naked eye. From this, we find that, while the variance is significantly reduced (reduction factors $11$, $12$ and $1219$), the variance-reduced estimator is biased for the estimated and averaged initializations. Moreover, the resulting estimator is only unbiased when using an exact initialization, leading to an even stronger reduction in variance.

\paragraph{HMM with variance reduction with Metropolis-Hastings} Since the exact initialization used in the previous experiment is generally not possible, we apply the Metropolis correction in the initialization to avoid the bias. The results are shown in figure \ref{fig:pdfs_nonlinear} (right), in which the blue and green dot-dashed lines correspond to estimator distributions using an estimated and averaged initialization, respectively. In this case, we derive the same conclusion as for the linear system (see figure \ref{fig:pdfs_linear}, right): the blue distribution is unbiased and shows a clear reduction in variance, while the green distribution is biased and gives only little reduction. To conclude, we investigate the estimated and averaged initializations with Metropolis-Hastings in figure \ref{fig:pdfs_nonlinear_est_av_MH}. It is seen that using $\FHMM{X^0}{500}{\omega_{0,j}}$ and $\FHMM{X^0}{5000}{\omega_{0,j}}$ lead to an unbiased variance-reduced estimator with reduction factors $10$ and $104$ (left plot), while using $\FHMM{X^0}{S,M}{\omega_{0,j}}$ becomes unbiased for fixed $M$ and $S \to \infty$ with reduction factors $1$, $2$ and $22$ (right plot).

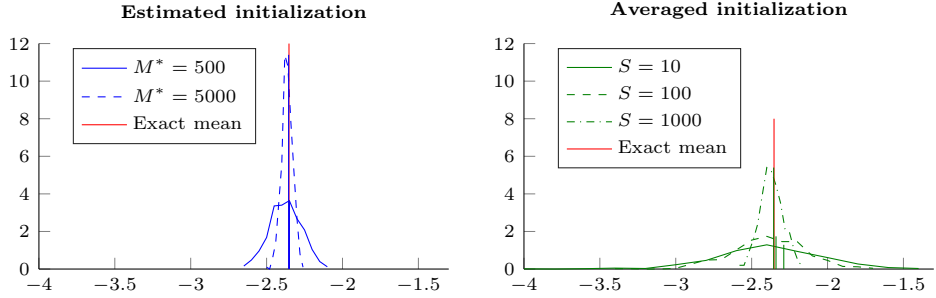
\begin{figure}[t]
	\begin{center}
		\figname{pdfs_nonlinear_est_av_MH}
%
%
\definecolor{mycolor1}{rgb}{0.00000,0.49804,0.00000}%
\newcommand{\figWidth}{0.33\textwidth} 
\newcommand{\figHeight}{3cm} 
\newcommand{\figSpacingRight}{1cm} 
\begin{tikzpicture}

\begin{axis}[%
width=0,
height=0,
scale only axis,
clip=false,
xmin=0,
xmax=1,
ymin=0,
ymax=1,
hide axis,
]
\node[align=center, text=black]
at (\figWidth+\figSpacingRight/2,\figHeight+0.9cm) {Variance-reduced estimator distributions with MH (nonlinear case)};
\end{axis}

\begin{axis}[%
width=\figWidth,
height=\figHeight,
at={(0,0.0)},
scale only axis,
xmin=-4,
xmax=-1.3,
ymin=0,
ymax=12,
ytick={0, 2, ..., 12},
axis background/.style={fill=white},
title style={font=\scriptsize\bfseries},
title={Estimated initialization},
axis x line*=bottom,
axis y line*=left,
legend style={shift={(-0.45,0)},legend cell align=left,align=left,draw=white!15!black}
]
\addplot [color=blue,solid]
  table[]{tikz/data/nonlinear_est_av_MH_pdfs-8.tsv};
\addlegendentry{$M^* = 500$};

\addplot [color=blue,dashed]
  table[]{tikz/data/nonlinear_est_av_MH_pdfs-9.tsv};
\addlegendentry{$M^* = 5000$};

\addplot [color=red,solid]
  table[]{tikz/data/nonlinear_est_av_MH_pdfs-10.tsv};
\addlegendentry{Exact mean};

\addplot [color=blue,solid,forget plot]
  table[]{tikz/data/nonlinear_est_av_MH_pdfs-11.tsv};
\addplot [color=blue,solid,forget plot]
  table[]{tikz/data/nonlinear_est_av_MH_pdfs-12.tsv};
\end{axis}

\begin{axis}[%
width=\figWidth,
height=\figHeight,
at={(\figWidth+\figSpacingRight,0)},
scale only axis,
xmin=-4,
xmax=-1.3,
ymin=0,
ymax=12,
ytick={0, 2, ..., 12},
axis background/.style={fill=white},
title style={font=\scriptsize\bfseries},
title={Averaged initialization},
axis x line*=bottom,
axis y line*=left,
legend style={shift={(-0.45,0)},legend cell align=left,align=left,draw=white!15!black}
]
\addplot [color=mycolor1,solid]
  table[]{tikz/data/nonlinear_est_av_MH_pdfs-1.tsv};
\addlegendentry{$S = 10$};

\addplot [color=mycolor1,dashed]
  table[]{tikz/data/nonlinear_est_av_MH_pdfs-2.tsv};
\addlegendentry{$S = 100$};

\addplot [color=mycolor1,dashdotted]
  table[]{tikz/data/nonlinear_est_av_MH_pdfs-3.tsv};
\addlegendentry{$S = 1000$};

\addplot [color=red,solid]
  table[]{tikz/data/nonlinear_est_av_MH_pdfs-4.tsv};
\addlegendentry{Exact mean};

\addplot [color=mycolor1,solid,forget plot]
  table[]{tikz/data/nonlinear_est_av_MH_pdfs-5.tsv};
\addplot [color=mycolor1,solid,forget plot]
  table[]{tikz/data/nonlinear_est_av_MH_pdfs-6.tsv};
\addplot [color=mycolor1,solid,forget plot]
  table[]{tikz/data/nonlinear_est_av_MH_pdfs-7.tsv};
\end{axis}
\end{tikzpicture}%
	\end{center}
	\vspace{-0.4cm}\caption{\label{fig:pdfs_nonlinear_est_av_MH} Variance-reduced estimator distributions (including Metropolis-Hastings) after one time step $\Dt = 0.02$ when applied to the nonlinear system \eqref{eq:nonlin_slow_fast}. The exact mean given in \eqref{eq:bias_var_nonlin_mean_exact} is depicted by a vertical red line in each plot. Left: estimated initialization. Right: averaged initialization. }
\end{figure}

\subsection{Local variance reduction} \label{subsec:local_var_reduction}
In this section, we demonstrate the local reduction in variance of the different estimators by considering their variance after $N=1$, $N=4$ and $N=10$ time steps of the method. The variance is obtained by repeating the computations over 100 realizations. For a fixed number of time steps, we study the influence of the macroscopic time step $\Dt$ by varying it as:
\begin{equation}
	\Dt = [0.1, 0.05, 0.02, 0.01, 0.005, 0.002, 0.001].
\end{equation}
We apply the estimators both to the linear and nonlinear system using the same parameters as in sections \ref{subsubsec:bias_var_lin} and \ref{subsubsec:bias_var_nonlin}, respectively.

For the linear system \eqref{eq:lin_slow_fast}, the results are plotted in figure \ref{fig:ordertest_linear}. In each plot, the HMM estimator variance without variance reduction (blue dot-dashed line) is seen to be constant as a function of $\Dt$ and is unchanged when increasing the number of macroscopic time steps $N$. This is because As sthe variance of a Markov chain Monte Carlo method is essentially of the order $1/M$ (here, $M=50$), thus constant in $\Dt$ and independent of time. On the contrary, the variance-reduced estimator with $\FHMM{X^0}{500}{\omega_{0}}$ (solid blue) and $\FHMM{X^0}{10,50}{\omega_{0}}$ (solid green) as initial estimators yields a variance that clearly depends on $\Dt$ and $N$. As indicated on figure \ref{fig:pdfs_linear} (middle), both initializations lead to variances that are much alike. We observe that, for fixed $N$, the variance becomes smaller for increasing time step $\Dt$, while, for a fixed time step $\Dt$, the variance decays rapidly with increasing $N$. This behavior is confirmed by our analysis, see equation \eqref{eq:lin_F_var}, for which we have $\Btilde = 0$ due to our choice of $\dt$ in equation \eqref{eq:lin_nobias_condition}. In that case, since $\lambda + pq/A < 0$ due to stability, equation \eqref{eq:lin_F_var} shows that for fixed $N$ the variance $\Var\left[\FHMMbar{\Xbar^N}{M}\right]$ converges to $\Var\left[\FHMMbar{\Xbar^0}{M}\right]$ from below for $\Dt \to 0$. It also states that for fixed $\Dt$ the variance $\Var\left[\FHMMbar{\Xbar^N}{M}\right]$ decays exponentially with increasing $N$. The expected evolution of variance in equation \eqref{eq:lin_F_var} is depicted by a dashed red line in each plot.

\begin{figure}[t]
	\begin{center}
		\figname{ordertest_linear}
%
%
\definecolor{mycolor1}{rgb}{0.00000,0.49804,0.00000}%
\newcommand{\figWidth}{3.7cm} 
\newcommand{\figHeight}{3cm} 
\newcommand{\figSpacingRight}{1.5cm} 
\begin{tikzpicture}

\begin{axis}[%
width=0,
height=0,
scale only axis,
clip=false,
xmin=0,
xmax=1,
ymin=0,
ymax=1,
hide axis,
]
\node[align=center, text=black]
at (3*\figWidth/2+\figSpacingRight,\figHeight+0.9cm) {Local variance reduction (linear case)};
\end{axis}

\begin{axis}[%
width=\figWidth,
height=\figHeight,
at={(0,0.0)},
scale only axis,
xmode=log,
xmin=0.001,
xmax=0.1,
xminorticks=true,
xlabel={$\Dt$},
xmajorgrids,
xminorgrids,
ymode=log,
ymin=1e-16,
ymax=1,
ytick={1e-14,1e-10,1e-6,1e-2},
yminorticks=true,
ylabel={Variance},
ylabel style={yshift=-0.1cm},
ymajorgrids,
yminorgrids,
grid style={dotted},
axis background/.style={fill=white},
title style={font=\scriptsize\bfseries},
title={$N = 1$}
]
\addplot [color=blue,dashdotted,mark=o,mark options={solid},forget plot]
  table[]{tikz/data/linear_ordertest-1.tsv};
\addplot [color=blue,solid,mark=o,mark options={solid},forget plot]
  table[]{tikz/data/linear_ordertest-3.tsv};
\addplot [color=mycolor1,solid,mark=o,mark options={solid},forget plot]
  table[]{tikz/data/linear_ordertest-2.tsv};
\addplot [color=red,dashed,forget plot]
  table[]{tikz/data/linear_ordertest-4.tsv};
\end{axis}

\begin{axis}[%
width=\figWidth,
height=\figHeight,
at={(\figWidth+\figSpacingRight,0)},
scale only axis,
xmode=log,
xmin=0.001,
xmax=0.1,
xminorticks=true,
xlabel={$\Dt$},
xmajorgrids,
xminorgrids,
ymode=log,
ymin=1e-16,
ymax=1,
ytick={1e-14,1e-10,1e-6,1e-2},
yminorticks=true,
ylabel={Variance},
ylabel style={yshift=-0.1cm},
ymajorgrids,
yminorgrids,
grid style={dotted},
axis background/.style={fill=white},
title style={font=\scriptsize\bfseries},
title={$N = 4$}
]
\addplot [color=blue,dashdotted,mark=o,mark options={solid},forget plot]
  table[]{tikz/data/linear_ordertest-5.tsv};
\addplot [color=blue,solid,mark=o,mark options={solid},forget plot]
  table[]{tikz/data/linear_ordertest-7.tsv};
\addplot [color=mycolor1,solid,mark=o,mark options={solid},forget plot]
  table[]{tikz/data/linear_ordertest-6.tsv};
\addplot [color=red,dashed,forget plot]
  table[]{tikz/data/linear_ordertest-8.tsv};
\end{axis}

\begin{axis}[%
width=\figWidth,
height=\figHeight,
at={(2*\figWidth+2*\figSpacingRight,0)},
scale only axis,
xmode=log,
xmin=0.001,
xmax=0.1,
xminorticks=true,
xlabel={$\Dt$},
xmajorgrids,
xminorgrids,
ymode=log,
ymin=1e-16,
ymax=1,
ytick={1e-14,1e-10,1e-6,1e-2},
yminorticks=true,
ylabel={Variance},
ylabel style={yshift=-0.1cm},
ymajorgrids,
yminorgrids,
grid style={dotted},
axis background/.style={fill=white},
title style={font=\scriptsize\bfseries},
title={$N = 10$}
]
\addplot [color=blue,dashdotted,mark=o,mark options={solid},forget plot]
  table[]{tikz/data/linear_ordertest-9.tsv};
\addplot [color=blue,solid,mark=o,mark options={solid},forget plot]
  table[]{tikz/data/linear_ordertest-11.tsv};
\addplot [color=mycolor1,solid,mark=o,mark options={solid},forget plot]
  table[]{tikz/data/linear_ordertest-10.tsv};
\addplot [color=red,dashed,forget plot]
  table[]{tikz/data/linear_ordertest-12.tsv};
\end{axis}
\end{tikzpicture}%
	\end{center}
	\vspace{-0.4cm}\caption{\label{fig:ordertest_linear} Local variance reduction of estimators evaluated after one (left), four (middle) and ten (right) time steps as a function of the macroscopic time step $\Dt$ for linear system \eqref{eq:lin_slow_fast}. Blue dot-dashed line: HMM estimator without variance reduction; solid blue and green lines: variance-reduced HMM estimator using an estimated and averaged initialization, respectively; red dashed line: expected variance according to \eqref{eq:lin_F_var}. }
\end{figure}
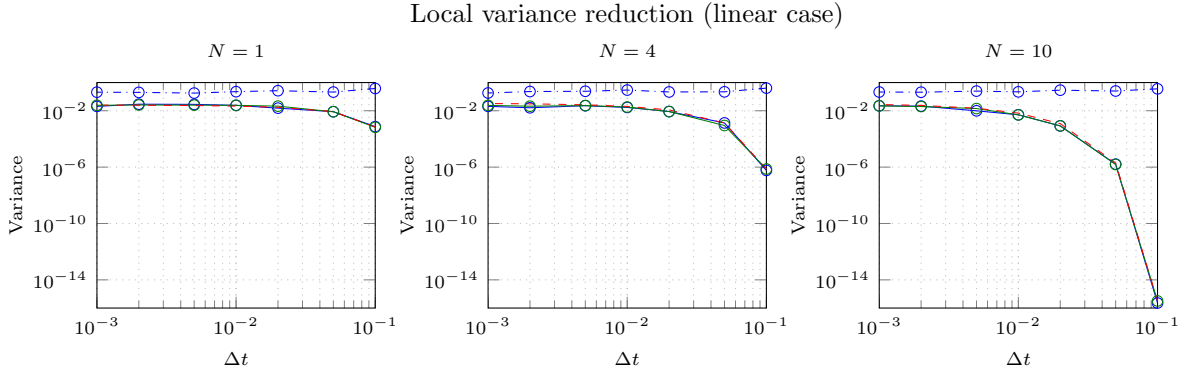

When repeating the above experiment for the nonlinear system \eqref{eq:nonlin_slow_fast}, we obtain the plots in figure \ref{fig:ordertest_nonlinear}. Once more, the HMM estimator variance (blue dot-dashed line) is roughly constant in $\Dt$ and $N$. The red, blue and green solid lines correspond to the variance of the variance-reduced estimator when using $F(X^0)$, $\FHMM{X^0}{500}{\omega_{0}}$ and $\FHMM{X^0}{10,50}{\omega_{0}}$ as initialization, respectively. The red line represents the best possible performance of the proposed variance reduction method. It confirms the formal result obtained in equation \eqref{eq:F_cv_var}, that is: (i) the variance increases with increasing $N$; (ii) although $\Var[\FHMMbar{X^0}{M}] = 0$ there is a small contribution to the variance due to the dependence of $A_N^j$ and $B_N^j$ on the Brownian path which is observed for $N=1$; (iii) since a little bit of variance is added in every macroscopic step, for $N=1$ the variance decays as $\Dt^2$ and this slope gradually decreases for increasing $N$. The solid blue and green line show that the estimated and averaged initializations give rise to a constant variance in $\Dt$ and $N$ corresponding to the variance of the initial estimator. As observed in figure \ref{fig:ordertest_nonlinear} (right), the averaged initialization resulted in practically no reduction in variance, which is clearly visible in each plot of figure \ref{fig:ordertest_nonlinear}.

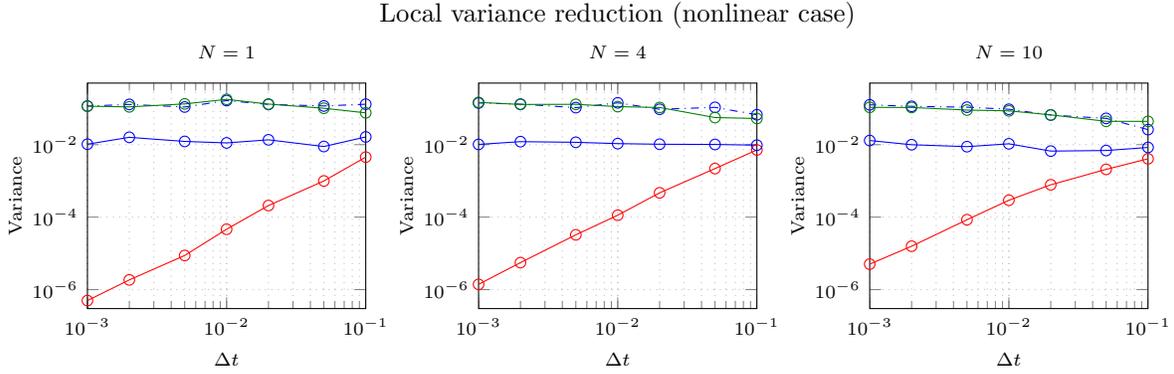
\begin{figure}[t]
	\begin{center}
		\figname{ordertest_nonlinear}
%
%
\definecolor{mycolor1}{rgb}{0.00000,0.49804,0.00000}%
\newcommand{\figWidth}{3.7cm} 
\newcommand{\figHeight}{3cm} 
\newcommand{\figSpacingRight}{1.5cm} 
\begin{tikzpicture}

\begin{axis}[%
width=0,
height=0,
scale only axis,
clip=false,
xmin=0,
xmax=1,
ymin=0,
ymax=1,
hide axis,
]
\node[align=center, text=black]
at (3*\figWidth/2+\figSpacingRight,\figHeight+0.9cm) {Local variance reduction (nonlinear case)};
\end{axis}

\begin{axis}[%
width=\figWidth,
height=\figHeight,
at={(0,0.0)},
scale only axis,
xmode=log,
xmin=0.001,
xmax=0.1,
xminorticks=true,
xlabel={$\Dt$},
xmajorgrids,
xminorgrids,
ymode=log,
ymin=3e-07,
ymax=0.5,
yminorticks=true,
ylabel={Variance},
ylabel style={yshift=-0.1cm},
ymajorgrids,
yminorgrids,
grid style={dotted},
axis background/.style={fill=white},
title style={font=\scriptsize\bfseries},
title={$N = 1$}
]
\addplot [color=red,solid,mark=o,mark options={solid},forget plot]
  table[]{tikz/data/nonlinear_ordertest-1.tsv};
\addplot [color=blue,solid,mark=o,mark options={solid},forget plot]
  table[]{tikz/data/nonlinear_ordertest-2.tsv};
\addplot [color=blue,dashdotted,mark=o,mark options={solid},forget plot]
  table[]{tikz/data/nonlinear_ordertest-3.tsv};
\addplot [color=mycolor1,solid,mark=o,mark options={solid},forget plot]
  table[]{tikz/data/nonlinear_ordertest-4.tsv};
\end{axis}

\begin{axis}[%
width=\figWidth,
height=\figHeight,
at={(\figWidth+\figSpacingRight,0)},
scale only axis,
xmode=log,
xmin=0.001,
xmax=0.1,
xminorticks=true,
xlabel={$\Dt$},
xmajorgrids,
xminorgrids,
ymode=log,
ymin=3e-07,
ymax=0.5,
yminorticks=true,
ylabel={Variance},
ylabel style={yshift=-0.1cm},
ymajorgrids,
yminorgrids,
grid style={dotted},
axis background/.style={fill=white},
title style={font=\scriptsize\bfseries},
title={$N = 4$}
]
\addplot [color=blue,dashdotted,mark=o,mark options={solid},forget plot]
  table[]{tikz/data/nonlinear_ordertest-5.tsv};
\addplot [color=red,solid,mark=o,mark options={solid},forget plot]
  table[]{tikz/data/nonlinear_ordertest-6.tsv};
\addplot [color=blue,solid,mark=o,mark options={solid},forget plot]
  table[]{tikz/data/nonlinear_ordertest-7.tsv};
\addplot [color=mycolor1,solid,mark=o,mark options={solid},forget plot]
  table[]{tikz/data/nonlinear_ordertest-8.tsv};
\end{axis}

\begin{axis}[%
width=\figWidth,
height=\figHeight,
at={(2*\figWidth+2*\figSpacingRight,0)},
scale only axis,
xmode=log,
xmin=0.001,
xmax=0.1,
xminorticks=true,
xlabel={$\Dt$},
xmajorgrids,
xminorgrids,
ymode=log,
ymin=3e-07,
ymax=0.5,
yminorticks=true,
ylabel={Variance},
ylabel style={yshift=-0.1cm},
ymajorgrids,
yminorgrids,
grid style={dotted},
axis background/.style={fill=white},
title style={font=\scriptsize\bfseries},
title={$N = 10$}
]
\addplot [color=blue,dashdotted,mark=o,mark options={solid},forget plot]
  table[]{tikz/data/nonlinear_ordertest-9.tsv};
\addplot [color=red,solid,mark=o,mark options={solid},forget plot]
  table[]{tikz/data/nonlinear_ordertest-10.tsv};
\addplot [color=blue,solid,mark=o,mark options={solid},forget plot]
  table[]{tikz/data/nonlinear_ordertest-11.tsv};
\addplot [color=mycolor1,solid,mark=o,mark options={solid},forget plot]
  table[]{tikz/data/nonlinear_ordertest-12.tsv};
\end{axis}
\end{tikzpicture}%
	\end{center}
	\vspace{-0.4cm}\caption{\label{fig:ordertest_nonlinear} Local variance reduction of estimators evaluated after one (left), four (middle) and ten (right) time steps as a function of the macroscopic time step $\Dt$ for nonlinear system \eqref{eq:nonlin_slow_fast}. Blue dot-dashed line: HMM estimator without variance reduction; solid red, blue and green lines: variance-reduced HMM estimator using an exact, estimated and averaged initialization, respectively. }
\end{figure}

\subsection{Solution trajectories} \label{subsec:trajectories}
Finally, we look at the solution paths of the averaged equation \eqref{eq:averaged_eq} obtained by different estimators for $F(X)$. We apply the method both to the linear (section \ref{subsubsec:trajectories_linear}) and nonlinear (section \ref{subsubsec:trajectories_nonlinear}) system.

\subsubsection{Linear system} \label{subsubsec:trajectories_linear}
Here, we focus on approximating the reduced evolution of the slow variable in equation \eqref{eq:lin_averaged_eq}, which is, in turn, an approximation of the slow variable's true evolution described in system \eqref{eq:lin_slow_fast}.
We compute the solution for $t \in [0,1]$ using initial conditions $X^0 = x_0 = 1$ and $y_0 = 1$. The system parameters in equation \eqref{eq:lin_slow_fast} are as follows: $\lambda = -10$, $p = 4$, $q = 0.5$ and $A = 1.2$.

We begin by applying the HMM technique without variance reduction generating $M=50$ samples in each iteration. For stability, the time step $\dt$ used in the Euler-Maruyama discretization of the fast equation is chosen as $\dt = \epsi$ with $\epsi=10^{-3}$. The forward Euler time step used in the discretization of the macroscopic equation \eqref{eq:averaged_eq} is fixed as $\Dt=0.02$. The time evolution of the variables $X$ and $F$ is depicted by the blue line in the left plots of figure \ref{fig:lin_paths}. The red line represents the exact solution of the macroscopic equation given in equation \eqref{eq:lin_X_exact}. Clearly, the statistical error dominates, thus justifying the need for variance reduction. The variance on $X$ and $F$ can be seen by the blue line in the right plots of figure \ref{fig:lin_paths} and is calculated by repeating the above experiment 100 times. We observe that the variance of the HMM estimator remains constant in time and behaves as $O(1/M)$ which is typical for a Markov chain Monte Carlo estimator.

\begin{figure}[t]
	\begin{center}
		\figname{paths_linear}
%
%
\definecolor{mycolor1}{rgb}{0.00000,0.49804,0.00000}%
\definecolor{mycyan}{rgb}{0.00000,0.74902,0.74902}%
\newcommand{\figWidth}{6cm} 
\newcommand{\figHeight}{2.5cm} 
\newcommand{\figSpacingRight}{2cm} 
\newcommand{\figSpacingTop}{1cm}
\begin{tikzpicture}

\begin{axis}[%
width=0,
height=0,
scale only axis,
clip=false,
xmin=0,
xmax=1,
ymin=0,
ymax=1,
hide axis,
]
\node[align=center, text=black]
at (\figWidth+1/2*\figSpacingRight,2*\figHeight+\figSpacingTop+0.5cm) {Solution trajectories and variance (linear case)};
\end{axis}

\begin{axis}[%
width=\figWidth,
height=\figHeight,
at={(0,\figHeight+\figSpacingTop)},
scale only axis,
xmin=0,
xmax=1,
xlabel={time},
ymin=-0.1,
ymax=1.1,
ytick={0, 0.2, ..., 1},
ylabel={$X(t)$},
axis background/.style={fill=white},
legend style={shift={(-0.02,-0.05)},legend cell align=left,align=left,draw=white!15!black}
]
\addplot [color=blue,solid,mark=*,mark size=1,mark options={solid}]
  table[]{tikz/data/lin_paths-10.tsv};
\addlegendentry{HMM};

\addplot [color=mycolor1,solid,mark=*,mark size=1,mark options={solid}]
  table[]{tikz/data/lin_paths-11.tsv};
\addlegendentry{VR HMM $(M^* = 500)$};

\addplot [color=red,solid,line width=1.0pt]
  table[]{tikz/data/lin_paths-12.tsv};
\addlegendentry{Exact};

\addplot [dashed] table[]{tikz/data/lin_paths-16.tsv};
\end{axis}

\begin{axis}[%
width=\figWidth,
height=\figHeight,
at={(0,0)},
scale only axis,
xmin=0,
xmax=1,
xlabel={time},
ymin=-8.5,
ymax=1.5,
ylabel={$F(X)$},
axis background/.style={fill=white}
]
\addplot [color=blue,solid,mark=*,mark size=1,mark options={solid},forget plot]
  table[]{tikz/data/lin_paths-13.tsv};
\addplot [color=mycolor1,solid,mark=*,mark size=1,mark options={solid},forget plot]
  table[]{tikz/data/lin_paths-14.tsv};
\addplot [color=red,solid,line width=1.0pt,forget plot]
  table[]{tikz/data/lin_paths-15.tsv};
\addplot [dashed] table[]{tikz/data/lin_paths-16.tsv};
\end{axis}

\begin{axis}[%
width=\figWidth,
height=\figHeight/2,
at={(\figWidth+\figSpacingRight,\figHeight+\figSpacingTop)},
scale only axis,
xmin=0,
xmax=1,
xmajorgrids,
xminorgrids,
xlabel={time},
ymode=log,
ymin=1e-32,
ymax=1e-25,
ytick={1e-31,1e-28},
ymajorgrids,
yminorgrids,
grid style={dotted},
axis background/.style={fill=white},
axis x line*=bottom,
]
\addplot [color=red,solid,mark=*,mark size=1,mark options={solid},forget plot]
  table[]{tikz/data/lin_paths-1.tsv};
\end{axis}

\begin{axis}[%
width=\figWidth,
height=\figHeight/2,
at={(\figWidth+\figSpacingRight,3/2*\figHeight+\figSpacingTop)},
scale only axis,
xmin=0,
xmax=1,
every x tick label/.append style={font=\color{white}},
xmajorgrids,
xminorgrids,
ymode=log,
ymin=1e-06,
ymax=0.004,
ytick={1e-6,1e-5,1e-4,1e-3},
ymajorgrids,
ylabel={$\Var[X]$},
ylabel style={xshift=-1cm},
grid style={dotted},
axis background/.style={fill=white},
axis x line*=top,
legend style={shift={(-0.05,1)},legend cell align=left,align=left,draw=white!15!black},
legend columns=2,
transpose legend
]
\addplot [color=blue,solid,mark=*,mark size=1,mark options={solid}]
  table[]{tikz/data/lin_paths-3.tsv};
\addlegendentry{HMM};
\addplot [color=mycolor1,solid,mark=*,mark size=1,mark options={solid}]
  table[]{tikz/data/lin_paths-4.tsv};
\addlegendentry{$M^* = 500$};
\addplot [color=mycyan,solid,mark=*,mark size=1,mark options={solid}]
  table[]{tikz/data/lin_paths-5.tsv};
\addlegendentry{$M^* = 5000$};
\addplot [color=red,solid,mark=*,mark size=1,mark options={solid}]
  table[]{tikz/data/lin_paths-6.tsv};
\addlegendentry{Exact init};
\end{axis}

\draw (\figWidth+\figSpacingRight-0.2cm,3/2*\figHeight+\figSpacingTop-0.2cm) -- ++(5:0.4cm);
\draw (\figWidth+\figSpacingRight-0.2cm,3/2*\figHeight+\figSpacingTop-0.3cm) -- ++(5:0.4cm);

\draw (2*\figWidth+\figSpacingRight-0.2cm,3/2*\figHeight+\figSpacingTop-0.2cm) -- ++(5:0.4cm);
\draw (2*\figWidth+\figSpacingRight-0.2cm,3/2*\figHeight+\figSpacingTop-0.3cm) -- ++(5:0.4cm);

\begin{axis}[%
width=\figWidth,
height=\figHeight/2,
at={(\figWidth+\figSpacingRight,0)},
scale only axis,
xmin=0,
xmax=1,
xmajorgrids,
xminorgrids,
xlabel={time},
ymode=log,
ymin=1e-29,
ymax=1e-23,
ytick={1e-28,1e-26},
ymajorgrids,
yminorgrids,
grid style={dotted},
axis background/.style={fill=white},
axis x line*=bottom,
]
\addplot [color=red,solid,mark=*,mark size=1,mark options={solid},forget plot]
  table[]{tikz/data/lin_paths-2.tsv};
\end{axis}

\begin{axis}[%
width=\figWidth,
height=\figHeight/2,
at={(\figWidth+\figSpacingRight,\figHeight/2)},
scale only axis,
xmin=0,
xmax=1,
every x tick label/.append style={font=\color{white}},
xmajorgrids,
xminorgrids,
ymode=log,
ymin=1e-12,
ymax=1,
ytick={1e-11,1e-6,1e-1},
ylabel={$\Var[F(X)]$},
ylabel style={xshift=-1cm},
ymajorgrids,
yminorgrids,
grid style={dotted},
axis background/.style={fill=white},
axis x line*=top
]
\addplot [color=blue,solid,mark=*,mark size=1,mark options={solid},forget plot]
  table[]{tikz/data/lin_paths-7.tsv};
\addplot [color=mycolor1,solid,mark=*,mark size=1,mark options={solid},forget plot]
  table[]{tikz/data/lin_paths-8.tsv};
\addplot [color=mycyan,solid,mark=*,mark size=1,mark options={solid}]
  table[]{tikz/data/lin_paths-9.tsv};
\end{axis}

\draw (\figWidth+\figSpacingRight-0.2cm,1/2*\figHeight-0.2cm) -- ++(5:0.4cm);
\draw (\figWidth+\figSpacingRight-0.2cm,1/2*\figHeight-0.3cm) -- ++(5:0.4cm);

\draw (2*\figWidth+\figSpacingRight-0.2cm,1/2*\figHeight-0.2cm) -- ++(5:0.4cm);
\draw (2*\figWidth+\figSpacingRight-0.2cm,1/2*\figHeight-0.3cm) -- ++(5:0.4cm);
\end{tikzpicture}%
	\end{center}
	\vspace{-0.4cm}\caption{\label{fig:lin_paths} Left: time evolution of $X$ and $F$ when applying HMM to the linear model problem \eqref{eq:lin_slow_fast} with and without variance reduction (green and blue lines, respectively). The red line represents the exact solution. Right: variance on $X$ and $F$. Blue line: HMM without variance reduction; red line: variance-reduced HMM with exact initialization; green and cyan lines: variance-reduced HMM with estimated initialization using $M^*=500$ and $M^*=5000$ samples, respectively. }
\end{figure}
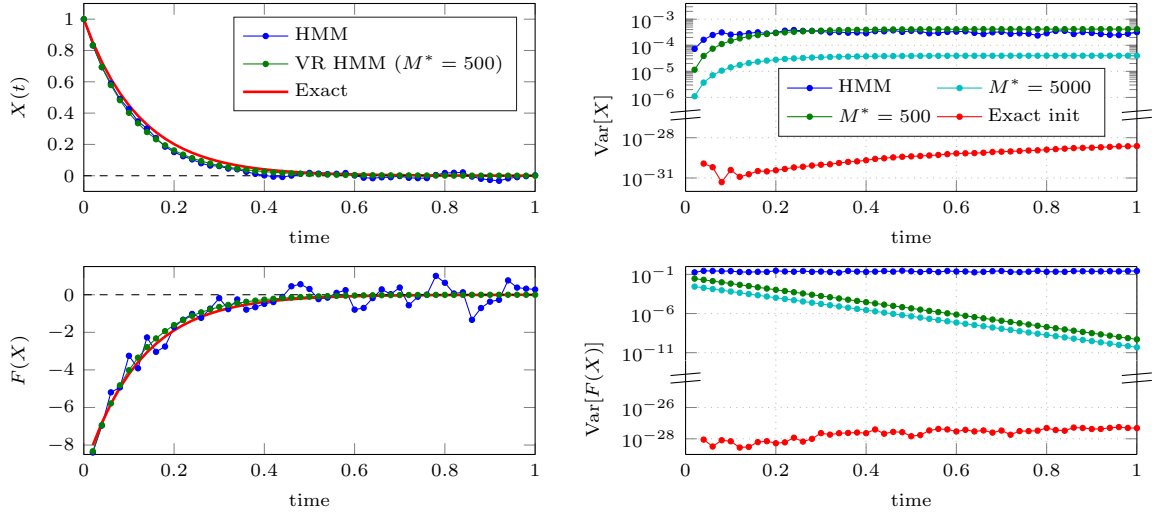

Next, we examine the proposed variance reduction technique based on control variables, as introduced in section \ref{subsec:hmm_cv}, for which we first need to specify the initial estimation $\FHMMbar{X^0}{M}$. As pointed out at the end of section \ref{subsec:est_var} and observed in section \ref{subsubsec:bias_var_lin}, when using an exact initialization $\FHMMbar{X^0}{M} = F(X^0)$ with $F(x)$ calculated in equation \eqref{eq:lin_averaged_eq}, the variance-reduced estimator is completely variance-free, thus leading to a deterministic estimator. This is confirmed by the right plots of figure \ref{fig:lin_paths}, which demonstrate that the variance on both $X$ and $F$ is indeed zero up to machine precision.
Subsequently, we compare the effect when using an estimated initialization $\FHMM{X^0}{M^*}{\omega_0}$. We remark that, to avoid introducing a bias in the linear case, the initial estimator is required to use a time step $\dt$ in the Euler-Maruyama scheme as derived in equation \eqref{eq:lin_nobias_condition}.
In all plots of figure \ref{fig:lin_paths}, the green line represents simulations when using $M^* = 500$ samples in the initial estimator. The left plots show that both $X$ and $F$ evolve much smoother than its HMM counterpart. The bottom right plot confirms that the variance on $F$ decays exponentially with time starting from the variance of the initial estimation, which was derived in equation \eqref{eq:lin_F_var}. However, the top right plot indicates that there is no reduction in variance on the trajectories $X$ with this initialization. To that end, when choosing $M^* = 5000$ samples in the initial estimator, we can improve the reduction in variance for both $X$ and $F$ by a factor $10$, which is shown by the cyan line on both right plots.

\subsubsection{Nonlinear system} \label{subsubsec:trajectories_nonlinear}
As a second model problem, we consider the nonlinear stochastic multiscale system given in equation \eqref{eq:nonlin_slow_fast}.
We calculate the solution of the averaged equation for $t \in [0,2]$ using initial conditions $X^0 = x_0 = 0.5$ and $y_0 = 0.5$.

We begin by applying the HMM procedure without variance reduction using $M=50$ samples by iterating over the Euler-Maruyama scheme \eqref{eq:EM_scheme} for the fast dynamics of system \eqref{eq:nonlin_slow_fast} with time step $\dt = \epsi$ and $\epsi = 10^{-3}$. We recall from section \ref{subsubsec:bias_var_nonlin} that the HMM solution converges to the wrong solution (that is, the solution of the wrong equation) in the nonlinear case, due to the time discretization error of the Euler-Maruyama scheme. Therefore, we require the Metropolis correction to avoid this bias. The forward Euler time step used in the discretization of the macroscopic equation is fixed as $\Dt = 0.05$. The time evolution of the variables $X$ and $F$ and their corresponding variance is depicted by the blue lines in figure \ref{fig:nonlin_paths}. The red line in the left plots represents the exact solution for $X$ and $F$ as given in \eqref{eq:nonlin_X_exact} and \eqref{eq:nonlin_F_exact}, respectively.

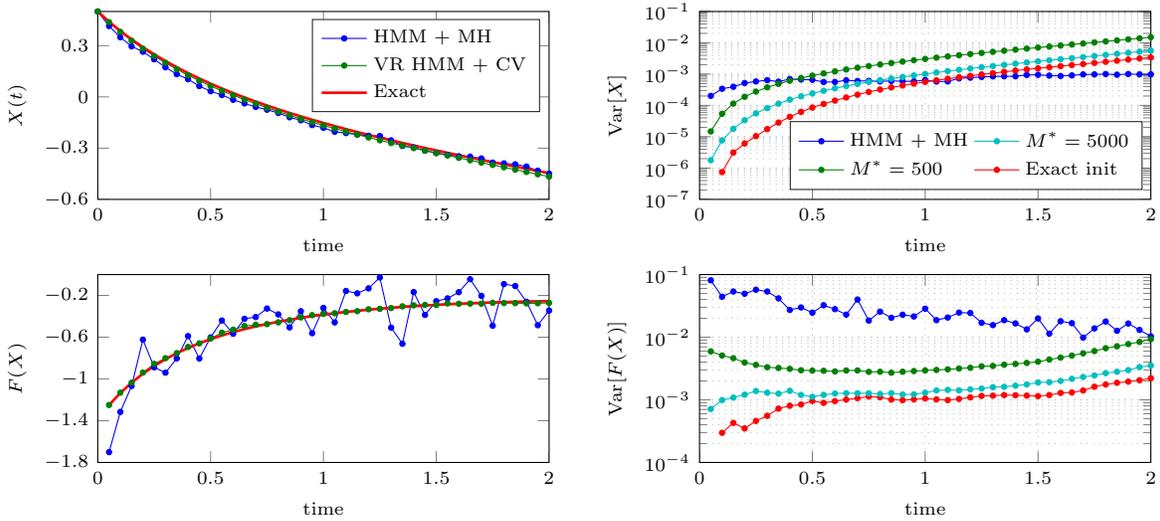
\begin{figure}[t]
	\begin{center}
		\figname{paths_nonlinear}
%
%
\definecolor{mygreen}{rgb}{0.00000,0.49804,0.00000}%
\definecolor{mypurple}{rgb}{0.49412,0.18431,0.55686}%
\definecolor{mycyan}{rgb}{0.00000,0.74902,0.74902}%
\newcommand{\figWidth}{6cm} 
\newcommand{\figHeight}{2.5cm} 
\newcommand{\figSpacingRight}{2cm} 
\newcommand{\figSpacingTop}{1cm}
\begin{tikzpicture}

\begin{axis}[%
width=0,
height=0,
scale only axis,
clip=false,
xmin=0,
xmax=1,
ymin=0,
ymax=1,
hide axis,
]
\node[align=center, text=black]
at (\figWidth+1/2*\figSpacingRight,2*\figHeight+\figSpacingTop+0.5cm) {Solution trajectories and variance (nonlinear case)};
\end{axis}

\begin{axis}[%
width=\figWidth,
height=\figHeight,
at={(0,\figHeight+\figSpacingTop)},
scale only axis,
xmin=0,
xmax=2,
xtick = {0, 0.5, ..., 2},
xlabel={time},
ymin=-0.6,
ymax=0.5,
ytick = {-0.6, -0.3, ..., 0.5},
yticklabel style={/pgf/number format/.cd,fixed},
ylabel={$X(t)$},
axis background/.style={fill=white},
legend style={shift={(0,-0.02)},legend cell align=left,align=left,draw=white!15!black}
]
\addplot [color=blue,solid,mark=*,mark size=1,mark options={solid}]
  table[]{tikz/data/nonlin_HMM_CV-1.tsv};
\addlegendentry{HMM + MH};

\addplot [color=mygreen,solid,mark=*,mark size=1,mark options={solid}]
  table[]{tikz/data/nonlin_HMM_CV-2.tsv};
\addlegendentry{VR HMM + CV};

\addplot [color=red,solid,line width=1pt]
  table[]{tikz/data/nonlin_HMM_CV-3.tsv};
\addlegendentry{Exact};

\end{axis}

\begin{axis}[%
width=\figWidth,
height=\figHeight,
at={(0,0)},
scale only axis,
xmin=0,
xmax=2,
xtick = {0, 0.5, ..., 2},
xlabel={time},
ymin=-1.8,
ymax=0,
ytick = {-1.8, -1.4, ..., 0},
yticklabel style={/pgf/number format/.cd,fixed},
ylabel={$F(X)$},
axis background/.style={fill=white}
]
\addplot [color=blue,solid,mark=*,mark size=1,mark options={solid},forget plot]
  table[]{tikz/data/nonlin_HMM_CV-6.tsv};
\addplot [color=mygreen,solid,mark=*,mark size=1,mark options={solid},forget plot]
  table[]{tikz/data/nonlin_HMM_CV-7.tsv};
\addplot [color=red,solid,line width=1pt,forget plot]
  table[]{tikz/data/nonlin_HMM_CV-8.tsv};
\end{axis}

\begin{axis}[%
width=\figWidth,
height=\figHeight,
at={(\figWidth+\figSpacingRight,\figHeight+\figSpacingTop)},
scale only axis,
xmin=0,
xmax=2,
xtick = {0, 0.5, ..., 2},
xlabel={time},
xmajorgrids,
ymode=log,
ymin=1e-07,
ymax=0.1,
yminorticks=true,
ytick = {1e-7, 1e-6, 1e-5, 1e-4, 1e-3, 1e-2, 1e-1},
ylabel={$\Var[X]$},
ymajorgrids,
yminorgrids,
grid style={dotted},
axis background/.style={fill=white},
legend style={at={(0.2,0.05)},anchor=south west,legend cell align=left,align=left,draw=white!15!black},
legend columns=2,
transpose legend
]
\addplot [color=blue,solid,mark=*,mark size=1,mark options={solid}]
  table[]{tikz/data/nonlin_HMM_CV-4.tsv};
\addlegendentry{HMM + MH};

\addplot [color=mygreen,solid,mark=*,mark size=1,mark options={solid}]
table[]{tikz/data/nonlin_HMM_CV2-2.tsv};
\addlegendentry{$M^*=500$};

\addplot [color=mycyan,solid,mark=*,mark size=1,mark options={solid}]
table[]{tikz/data/nonlin_HMM_CV2-3.tsv};
\addlegendentry{$M^*=5000$};

\addplot [color=red,solid,mark=*,mark size=1,mark options={solid}]
  table[]{tikz/data/nonlin_HMM_CV-5.tsv};
\addlegendentry{Exact init};

\end{axis}

\begin{axis}[%
width=\figWidth,
height=\figHeight,
at={(\figWidth+\figSpacingRight,0)},
scale only axis,
xmin=0,
xmax=2,
xtick = {0, 0.5, ..., 2},
xlabel={time},
xmajorgrids,
ymode=log,
ymin=0.0001,
ymax=0.1,
yminorticks=true,
ylabel={$\Var[F(X)]$},
ymajorgrids,
yminorgrids,
grid style={dotted},
axis background/.style={fill=white}
]
\addplot [color=blue,solid,mark=*,mark size=1,mark options={solid},forget plot]
  table[]{tikz/data/nonlin_HMM_CV-9.tsv};
\addplot [color=red,solid,mark=*,mark size=1,mark options={solid},forget plot]
  table[]{tikz/data/nonlin_HMM_CV-10.tsv};
\addplot [color=mygreen,solid,mark=*,mark size=1,mark options={solid},forget plot]
  table[]{tikz/data/nonlin_HMM_CV2-5.tsv};
\addplot [color=mycyan,solid,mark=*,mark size=1,mark options={solid},forget plot]
  table[]{tikz/data/nonlin_HMM_CV2-6.tsv};
\end{axis}
\end{tikzpicture}%
	\end{center}
	\vspace{-0.4cm}\caption{\label{fig:nonlin_paths} Left: time evolution of $X$ and $F$ when applying HMM with Metropolis-Hastings to the nonlinear model problem \eqref{eq:nonlin_slow_fast} with and without variance reduction (green and blue lines, respectively). The red line represents the exact solution. Right: variance on $X$ and $F$. Blue line: HMM without variance reduction; red line: variance-reduced HMM with exact initialization; green and cyan lines: variance-reduced HMM with estimated initialization using $M^* = 500$ and $M^* = 5000$ samples, respectively. }
\end{figure}

Next, we apply the variance-reduced HMM estimator. As noted in section \ref{subsubsec:est_bias_general}, when combining the Metropolis-Hastings algorithm with the variance-reduced estimator, we lose strong correlation between the HMM estimators $\FHMM{\Xhat^N}{M}{\omega_N}$ and $\FHMM{\Xhat^{N-1}}{M}{\omega_N}$ in equation \eqref{eq:F_cv}. This is clarified as follows. Both HMM estimators generate an ensemble of samples using the same seed $\omega_n$ but a different value of the slow variable. In general, it is not known a priori when and which samples will be accepted or rejected in the Metropolis-Hastings algorithm. Since samples can be rejected in different places in both ensembles, they are in principle no longer correlated. Consequently, when subtracting these two estimators the statistical error will be larger than that of the individual estimators. To resolve this, we instead use two classical HMM estimators without the Metropolis-Hastings extension in equation \eqref{eq:F_cv}, each producing a bias. However, since both contain the same bias, subtraction yields a result of order $\Dt$ (that is, the distance between $\Xbar^{N-1}$ and $\Xbar^N$) which lies within the accuracy of the forward Euler method.

In the variance-reduced setting, we perform experiments using both an exact initialization $F(X^0)$, as given in equation \eqref{eq:nonlin_F_exact}, as well as an estimated initialization $\FHMM{X^0}{M^*}{\omega_0}$ in the first forward Euler step of the averaged equation. The results are shown in figure \ref{fig:nonlin_paths}, where the green lines depict the evolution of the variables $X$ and $F$ and their corresponding variance when using $M^* = 500$ in the initial estimator. In addition, in the right plots we also show the variance on $X$ and $F$ when choosing an exact initialization (red lines) and when setting $M^* = 5000$ (cyan lines).
From this, we observe that there is a clear buildup of variance for all variance-reduced estimators.

As suggested at the end of section \ref{subsec:hmm_cv}, the variance buildup can be countered by occasionally reinitializing the estimator. Therefore, in what follows, we consider the influence of reinitializing the estimator after every $R$ macroscopic time steps; that is, we compute a new accurate estimation (similar to the initialization) after a fixed number of macroscopic time steps using $M_r$ samples. The resulting variance on $X$ and $F$ is plotted in figure \ref{fig:nonlin_paths_R}, where we compare the HMM estimator using $M = 50$ samples (solid blue line) with the variance-reduced estimator using $\FHMM{X^0}{500}{\omega_0}$ as initial estimator, $M = 20$ samples in the HMM estimator difference in \eqref{eq:F_cv}  for $R=20$ (green line), $R=10$ (red line), $R=5$ (cyan line) and $R=2$ (purple line) based on $M_r = 500$ samples. It is seen that the variance on the solution trajectories clearly depends on the value of $R$.
This experiment shows that, by repeatedly reinitializing the estimator, we can control the variance on both $X$ and $F$ and counter the buildup of variance on both quantities that was seen in figure \ref{fig:nonlin_paths}.

\begin{figure}[t]
	\begin{center}
		\figname{paths_nonlinear_R}
%
%
\definecolor{mygreen}{rgb}{0.00000,0.49804,0.00000}%
\definecolor{mycyan}{rgb}{0.00000,0.74902,0.74902}%
\definecolor{mypurple}{rgb}{1.00000,0.00000,1.00000}%
\newcommand{\figWidth}{6cm} 
\newcommand{\figHeight}{2.5cm} 
\newcommand{\figSpacingRight}{2cm} 
\newcommand{\figSpacingTop}{1cm}
\begin{tikzpicture}

\begin{axis}[%
width=0,
height=0,
scale only axis,
clip=false,
xmin=0,
xmax=1,
ymin=0,
ymax=1,
hide axis,
]
\node[align=center, text=black]
at (\figWidth+1/2*\figSpacingRight,2*\figHeight+\figSpacingTop+0.5cm) {Solution trajectories and variance with reinitialization (nonlinear case)};
\end{axis}

\begin{axis}[%
width=\figWidth,
height=\figHeight,
at={(0,\figHeight+\figSpacingTop)},
scale only axis,
xmin=0,
xmax=2,
xtick = {0, 0.5, ..., 2},
xlabel={time},
ymin=-0.6,
ymax=0.5,
ytick = {-0.6, -0.3, ..., 0.5},
yticklabel style={/pgf/number format/.cd,fixed},
ylabel={$X(t)$},
axis background/.style={fill=white},
legend style={shift={(-0.02,-0.05)},legend cell align=left,align=left,draw=white!15!black},
legend columns=2,
transpose legend
]
\addplot [color=mygreen,solid,mark=*,mark size=1,mark options={solid}]
  table[]{tikz/data/nonlin_paths_refresh-1.tsv};
\addlegendentry{$R = 20$};

\addplot [color=red,solid,mark=*,mark size=1,mark options={solid}]
  table[]{tikz/data/nonlin_paths_refresh-2.tsv};
\addlegendentry{$R = 10$};

\addplot [color=mycyan,solid,mark=*,mark size=1,mark options={solid}]
  table[]{tikz/data/nonlin_paths_refresh-3.tsv};
\addlegendentry{$R = 5$};

\addplot [color=black,solid,line width=1.0pt]
  table[]{tikz/data/nonlin_paths_R-4.tsv};
\addlegendentry{Exact};

\end{axis}

\begin{axis}[%
width=\figWidth,
height=\figHeight,
at={(0,0)},
scale only axis,
xmin=0,
xmax=2,
xtick = {0, 0.5, ..., 2},
xlabel={time},
ymin=-1.8,
ymax=0,
ytick = {-1.8, -1.4, ..., 0},
yticklabel style={/pgf/number format/.cd,fixed},
ylabel={$F(X)$},
axis background/.style={fill=white}
]
\addplot [color=black,solid,line width=1.0pt,forget plot]
  table[]{tikz/data/nonlin_paths_R-8.tsv};
\addplot [color=mygreen,solid,mark=*,mark size=1,mark options={solid},forget plot]
  table[]{tikz/data/nonlin_paths_refresh-4.tsv};
\addplot [color=red,solid,mark=*,mark size=1,mark options={solid},forget plot]
  table[]{tikz/data/nonlin_paths_refresh-5.tsv};
\addplot [color=mycyan,solid,mark=*,mark size=1,mark options={solid},forget plot]
  table[]{tikz/data/nonlin_paths_refresh-6.tsv};
\end{axis}

\begin{axis}[%
width=\figWidth,
height=\figHeight,
at={(\figWidth+\figSpacingRight,\figHeight+\figSpacingTop)},
scale only axis,
xmin=0,
xmax=2,
xtick = {0, 0.5, ..., 2},
xlabel={time},
xmajorgrids,
ymode=log,
ymin=1e-07,
ymax=0.1,
yminorticks=true,
ytick = {1e-7, 1e-6, 1e-5, 1e-4, 1e-3, 1e-2, 1e-1},
ylabel={$\Var[X]$},
ymajorgrids,
yminorgrids,
grid style={dotted},
axis background/.style={fill=white},
legend style={shift={(-0.02,0.6)},legend cell align=left,align=left,draw=white!15!black},
legend columns=2,
transpose legend
]
\addplot [color=blue,solid,mark=*,mark size=1,mark options={solid}]
  table[]{tikz/data/nonlin_HMM_CV-4.tsv};
\addlegendentry{HMM};
\addplot [color=mygreen,solid,mark=*,mark size=1,mark options={solid}]
  table[]{tikz/data/nonlin_paths_refresh-7.tsv};
\addlegendentry{$R = 20$};
\addplot [color=red,solid,mark=*,mark size=1,mark options={solid}]
  table[]{tikz/data/nonlin_paths_refresh-8.tsv};
\addlegendentry{$R = 10$};
\addplot [color=mycyan,solid,mark=*,mark size=1,mark options={solid}]
  table[]{tikz/data/nonlin_paths_refresh-9.tsv};
\addlegendentry{$R = 5$};
\addplot [color=mypurple,solid,mark=*,mark size=1,mark options={solid}]
  table[]{tikz/data/nonlin_paths_refresh_R2-1.tsv};
\addlegendentry{$R = 2$};
\end{axis}

\begin{axis}[%
width=\figWidth,
height=\figHeight,
at={(\figWidth+\figSpacingRight,0)},
scale only axis,
xmin=0,
xmax=2,
xtick = {0, 0.5, ..., 2},
xlabel={time},
xmajorgrids,
ymode=log,
ymin=0.0001,
ymax=0.1,
yminorticks=true,
ylabel={$\Var[F(X)]$},
ymajorgrids,
yminorgrids,
grid style={dotted},
axis background/.style={fill=white}
]
\addplot [color=blue,solid,mark=*,mark size=1,mark options={solid},forget plot]
  table[]{tikz/data/nonlin_HMM_CV-9.tsv};
\addplot [color=mygreen,solid,mark=*,mark size=1,mark options={solid},forget plot]
  table[]{tikz/data/nonlin_paths_refresh-10.tsv};
\addplot [color=red,solid,mark=*,mark size=1,mark options={solid},forget plot]
  table[]{tikz/data/nonlin_paths_refresh-11.tsv};
\addplot [color=mycyan,solid,mark=*,mark size=1,mark options={solid},forget plot]
  table[]{tikz/data/nonlin_paths_refresh-12.tsv};
\addplot [color=mypurple,solid,mark=*,mark size=1,mark options={solid},forget plot]
  table[]{tikz/data/nonlin_paths_refresh_R2-2.tsv};
\end{axis}
\end{tikzpicture}%
	\end{center}
	\vspace{-0.4cm}\caption{\label{fig:nonlin_paths_R} Left: time evolution of $X$ and $F$ when applying HMM with variance reduction to the nonlinear model problem \eqref{eq:nonlin_slow_fast} for different reinitialization frequencies $R$. The black line represents the exact solution. Right: variance on $X$ and $F$. Blue line: HMM without variance reduction; other lines: variance-reduced HMM using an estimated initialization with $M^* = 500$ and reinitializing after every $R$ time steps, for $R=20$ (green line), $R=10$ (red line), $R=5$ (cyan line) and $R=2$ (purple line). The reinitialization uses the same procedure as the initialization. }
\end{figure}
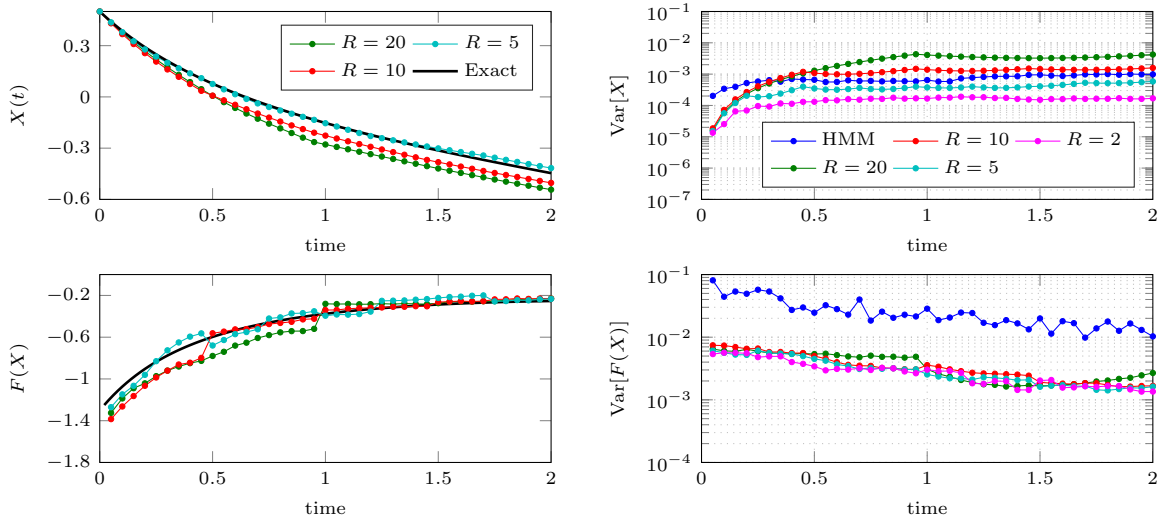

\section{Conclusions} \label{sec:conclusions}
We presented a variance reduction technique based on control variables for stochastic slow-fast systems containing a deterministic slow equation and a stochastic fast equation. The proposed method supplements the HMM estimator, which approximates the right hand side integral in the reduced description of the stochastic system, and we commented on its applicability for the coarse-projective integration estimator in the equation-free framework. 
We discussed the initialization of the variance-reduced estimator by means of an exact, estimated or averaged initial estimator and considered its reinitialization to counter the buildup of variance with time.
We analyzed the estimator variance and additional bias compared to the HMM estimator and derived explicit expressions for a linear stochastic system. We applied the variance-reduced estimator to a linear and nonlinear test problem, in which we considered the effects of the initialization and reinitialization and compared the results with the classical HMM estimator. The numerical experiments showed that nonlinear systems generally require the Metropolis-Hastings correction in both estimators to avoid introducing a bias. Moreover, reinitializing the estimator was found to be very effective to control the variance on solution trajectories. 

In the experiments, we did not compare in detail the computational cost of the different procedures for a desired variance. The variance-reduced scheme has many numerical parameters that can be chosen (such as the frequency and accuracy of reinitialization). Moreover, the standard HMM scheme benefits from a self-averaging effect, since the errors in subsequent macroscopic time steps are independent. In the variance-reduced scheme, the individual errors are smaller, but they are correlated.  As a consequence, a fair comparison of computational cost is highly non-trivial and may well be problem-dependent. We postpone such a comparison to future research.

\section*{References}

\bibliography{refs_mendeley}

\begin{thebibliography}{10}
\expandafter\ifx\csname url\endcsname\relax
  \def\url#1{\texttt{#1}}\fi
\expandafter\ifx\csname urlprefix\endcsname\relax\def\urlprefix{URL }\fi
\expandafter\ifx\csname href\endcsname\relax
  \def\href#1#2{#2} \def\path#1{#1}\fi

\bibitem{Berglund2006}
N.~Berglund, B.~Gentz, {Noise-induced phenomena in slow-fast dynamical systems:
  a sample-paths approach}, Springer Science {\&} Business Media, 2006.

\bibitem{Bruna2014}
M.~Bruna, S.~J. Chapman, M.~J. Smith, {Model reduction for slow-fast stochastic
  systems with metastable behaviour}, Journal of Chemical Physics 140~(17)
  (2014) 1--23.

\bibitem{Imkeller2001}
P.~Imkeller, J.-S. von Storch, {Stochastic Climate Models}, Birkh{\"{a}}user
  Basel, Basel, 2001.

\bibitem{Erban2009}
R.~Erban, S.~J. Chapman, I.~G. Kevrekidis, T.~Vejchodsk{\'{y}}, {Analysis of a
  stochastic chemical system close to a sniper bifurcation of its mean field
  model}, SIAM Journal on Applied Mathematics 70~(3) (2009) 984--1016.

\bibitem{Givon2004}
D.~Givon, R.~Kupferman, A.~Stuart, {Extracting macroscopic dynamics: model
  problems and algorithms}, Nonlinearity 17~(6) (2004) 55--127.

\bibitem{Rousset2013}
M.~Rousset, G.~Samaey, {Individual-Based Models for Bacterial Chemotaxis in the
  Diffusion Asymptotics}, Mathematical Models and Methods in Applied Sciences
  23~(11) (2013) 2005--2037.

\bibitem{Pavliotis2008}
G.~A. Pavliotis, A.~Stuart, {Multiscale methods: averaging and homogenization},
  Springer Science {\&} Business Media, 2008.

\bibitem{Li2008}
T.~Li, A.~Abdulle, W.~E, {Effectiveness of implicit methods for stiff
  stochastic differential equations}, Communications in Computational Physics
  3~(2) (2008) 295--307.

\bibitem{Vanden-Eijnden2003}
E.~Vanden-Eijnden, {Numerical techniques for multi-scale dynamical systems with
  stochastic effects}, Communications in Mathematical Sciences 1~(2) (2003)
  385--391.

\bibitem{E2005}
W.~E, D.~Liu, E.~Vanden-Eijnden, {Analysis of Multiscale Methods for Stochastic
  Differential Equations}, Communications on Pure and Applied Mathematics
  58~(11) (2005) 1544--1585.

\bibitem{caflisch1998monte}
R.~E. Caflisch, Monte carlo and quasi-monte carlo methods, Acta numerica 7
  (1998) 1--49.

\bibitem{E2003a}
W.~E, B.~Engquist, {The Heterogeneous Multiscale Methods}, Communications in
  Mathematical Sciences 1~(1) (2003) 87--132.

\bibitem{Abdulle2012}
A.~Abdulle, W.~E, B.~Engquist, E.~Vanden-Eijnden, {The heterogeneous multiscale
  method}, Acta Numerica 21~(May 2012) (2012) 1--87.

\bibitem{Givon2006}
D.~Givon, I.~G. Kevrekidis, R.~Kupferman, {Strong convergence of projective
  integration schemes for singularly perturbed stochastic differential
  systems}, Communications in Mathematical Sciences 4~(4) (2006) 707--729.

\bibitem{Kevrekidis2003}
I.~G. Kevrekidis, C.~W. Gear, J.~M. Hyman, P.~G. Kevrekidis, O.~Runborg,
  C.~Theodoropoulos, {Equation-Free, Coarse-Grained Multiscale Computation:
  enabling microscopic simulators to perform system-level tasks},
  Communications in Mathematical Sciences 1~(4) (2003) 715--762.

\bibitem{Kevrekidis2009}
I.~G. Kevrekidis, G.~Samaey, {Equation-free multiscale computation: algorithms
  and applications.}, Annual review of physical chemistry 60 (2009) 321--344.

\bibitem{Gear2002}
C.~W. Gear, I.~G. Kevrekidis, C.~Theodoropoulos, {'Coarse'
  integration/bifurcation analysis via microscopic simulators: micro-Galerkin
  methods}, Computers and Chemical Engineering 26~(7-8) (2002) 941--963.

\bibitem{melis2016variance}
W.~Melis, G.~Samaey, Variance-reduced {HMM} for stochastic slow-fast systems,
  Procedia Computer Science 80 (2016) 1255--1266.

\bibitem{Glasserman2003}
P.~Glasserman, {Monte Carlo Methods in Financial Engineering}, Springer Science
  {\&} Business Media, 2003.

\bibitem{Papavasiliou2007}
A.~Papavasiliou, I.~G. Kevrekidis, {Variance reduction for the equation-free
  simulation of multiscale stochastic systems}, Multiscale Modeling {\&}
  Simulation 6~(1) (2007) 70--89.

\bibitem{Higham2001}
D.~J. Higham, {An Algorithmic Introduction to Numerical Simulation of
  Stochastic Differential Equations}, SIAM Review 43~(3) (2001) 525--546.

\bibitem{roberts1996exponential}
G.~O. Roberts, R.~L. Tweedie, Exponential convergence of langevin distributions
  and their discrete approximations, Bernoulli (1996) 341--363.

\bibitem{Cances:EsaimM2An:2007}
E.~Canc\`{e}s, F.~Legoll, G.~Stoltz, Theoretical and numerical comparison of
  some sampling methods for molecular dynamics, ESAIM: M2AN 41~(2) (2007)
  351--389.

\bibitem{Fatkullin:JournalOfComputationalPhysics:2004}
I.~Fatkullin, E.~Vanden-Eijnden, A computational strategy for multiscale
  systems with applications to lorenz 96 model, Journal of Computational
  Physics 200~(2) (2004) 605--638.

\end{thebibliography}
	
\end{document}